\documentclass{article}
\usepackage{amssymb}

\newtheorem{definition}{Definition}[section]
\newtheorem{theorem}[definition]{Theorem}

\newtheorem{lemma}[definition]{Lemma}

\newtheorem{notation}[definition]{Notation}
\newtheorem{note}[definition]{Note}

\newtheorem{corollary}[definition]{Corollary}

\def\R{\mathbb R}
\def\C{\mathbb C}

\begin{document}
\title{The subconstituent algebra of  
a bipartite distance-regular \\ graph; thin   
modules with endpoint two}

\author{Mark MacLean\footnote{Mathematics Department, Seattle University, 901 Twelfth Avenue, Seattle WA  98122-1090 USA}  $\;$ and Paul Terwilliger\footnote{ Mathematics Department, University of Wisconsin, 480 Lincoln Drive, Madison WI  53706-1388 USA}
}

\date{}
\maketitle

\begin{abstract}
We consider a bipartite distance-regular graph $\Gamma $
with diameter 
$D \geq 4$, valency $k \geq 3$, intersection numbers $b_i, c_i$, distance matrices $A_i$, and  eigenvalues 
$\theta_0 > \theta_1 > \cdots > \theta_D$.
Let
$X$ denote the vertex set of $\Gamma $
and fix $x \in X$.
Let 
 $T=T(x)$
denote the subalgebra
of 
 $\hbox{Mat}_X(\C)$ 
generated by $A, E^*_0, E^*_1, \ldots, E^*_D$, where
$A=A_1$ and 
$E^*_i$ denotes the projection onto the $i^{\hbox{th}}$ 
subconstituent of $\Gamma $ with respect to $x$.
$T$ is called the {\it subconstituent algebra} (or {\it Terwilliger algebra}) of $\Gamma$ with respect to $x$.
An irreducible $T$-module $W$ is said to be {\it thin}
whenever $\hbox{dim}E^*_iW\leq 1$ for
$0 \leq i \leq D$.
By the {\it endpoint} of $W$ we mean 
$\hbox{min}\lbrace i |E^*_iW\not=0 \rbrace $.
Assume $W$ is thin with endpoint 2. Observe $E^*_2W$ is a 1-dimensional
eigenspace for $E^*_2A_2E^*_2$; let $\eta $ denote the
corresponding eigenvalue.  
It is known 
${\tilde \theta}_1  \leq 
\eta   
\leq {\tilde \theta}_d $ where 
 ${\tilde \theta}_1=-1-b_2 b_{3}(\theta_1^{2}-b_{2})^{-1},  \quad 
{\tilde \theta}_d=-1-b_2 b_{3}(\theta_d^{2}-b_{2})^{-1}$, and
$d = \lfloor D/2 \rfloor$.
To describe the structure of $W$ we distinguish four cases: (i) $\eta 
={\tilde \theta}_1 $; (ii) $D$ is odd 
and $\eta = {\tilde \theta}_{d}$; (iii) $D$ is even 
and $\eta = {\tilde \theta}_{d}$; (iv) $ {\tilde \theta}_1 < \eta < {\tilde 
\theta}_d$.  
We investigated cases (i), (ii) in \cite{maclean5}.  Here we investigate cases (iii), 
(iv) and obtain the following
results.  We show the dimension of $W$ is $D-1-e$ where
$e=1$ in case (iii) and $e=0$ in case (iv).  Let $v$ denote a nonzero vector in $E^*_2W$.
We show $W$ has a basis $E_iv$ $(i \in S)$,
where $E_i$ denotes the primitive idempotent of $A$ associated
with $\theta_i$ and where 
the set $S$ is $\{1, 2, \ldots, d-1\} \cup \{d+1, d+2, \ldots, D-1\}$ in case (iii)
and $\{1, 2, \ldots, D-1\}$ in case (iv).
We show this basis is orthogonal (with respect to the Hermitian dot
product) and 
we compute the square-norm of each basis vector. 
We show $W$ has a basis $E^*_{i+2}A_iv$  $(0 \leq i \leq D-2-e)$, and we  find 
the matrix representing $A$ with respect to this basis.
We show this basis is orthogonal 
and we compute the square-norm of each basis vector. 
We find the transition matrix relating our two bases for $W$.
\\ \\
{\bf Keywords}. Distance-regular graph, association
scheme, Terwilliger algebra, subconstituent algebra. \hfil\break
{\bf 2000 Mathematics Subject Classification}. Primary 05E30; Secondary 05E35, 05C50
\end{abstract}

\section{Introduction}

\bigskip \noindent
Let $\Gamma$ denote a distance-regular graph with diameter $D \geq 
4$, valency $k \ge 3$, intersection numbers $a_i, b_i, c_i$, and distance matrices $A_i$
(see Section 2 for formal definitions).
We recall the subconstituent algebra of $\Gamma$.
Let
$X$ denote the vertex set of $\Gamma$ and fix $x 
\in X$. We view $x$ as a ``base vertex.'' Let $T=T(x)$ denote the subalgebra of
$\hbox{Mat}_X(\C)$ generated by $A, E^*_0, E^*_1, \ldots, E^*_D$, 
where
$A=A_1$ and $E^*_i$ represents the 
projection onto the $i^{\hbox{th}}$ subconstituent of $\Gamma $ with 
respect to $x$. The algebra $T$ is called the {\it subconstituent 
algebra} (or {\it Terwilliger algebra}) of $\Gamma $
with respect to $x$
\cite{terwSub1}. Observe $T$ has finite dimension.
Moreover $T$ is semi-simple; the reason is
each of $A, E^*_0, E^*_1, \ldots, E^*_D$ is symmetric with
real entries, so
$T$ is closed under the conjugate-transpose map \cite[p. 157]{CR}. 
Since $T$ is semi-simple, each $T$-module is a direct sum of 
irreducible $T$-modules. Describing the irreducible $T$-modules is an
active area of research
\cite{caugh2}--\cite{curtin6}, \cite{dickie1}--\cite{go2}, 
\cite{hobart},
\cite{maclean5}--\cite{tomiyama2}.

\bigskip \noindent
In this paper we are concerned with the irreducible $T$-modules that 
possess a certain property.  In order to define this property we make 
a few observations.
Let $W$ 
denote an irreducible $T$-module.
Then $W$ is the direct sum of the
nonzero spaces among $E^*_0W, E^*_1W,\ldots, E^*_DW$. There is a 
second decomposition of interest.
To obtain it we make a definition.
Let $k=\theta_0 > \theta_1 > \cdots > \theta_D$
denote the distinct eigenvalues of $A$, and for $0 \leq i \leq D $ 
let $E_i$
denote the primitive idempotent of $A$ associated with $\theta_i$.
Then $W$ is the direct sum of the
nonzero spaces among $E_0W, E_1W,\ldots, E_DW$.
If the dimension of $E^*_iW$ is at most 1 for
$0 \leq i \leq D$ then
the dimension of $E_iW$ is at most 1 for
$0 \leq i \leq D$ \cite[Lemma 3.9]{terwSub1};
in this case we say $W$ is {\it thin}.
Let $W$ denote an irreducible $T$-module.
By the {\it endpoint} of $W$
we mean $\hbox{min}\lbrace i | 0 \leq i \leq D, \;E^*_iW\not=0 
\rbrace $.
There exists a unique irreducible $T$-module with endpoint 0 
\cite[Proposition 8.4]{egge1}. We call this module $V_0$.
The module $V_0$ is thin; in fact
$E^*_iV_0$ and
$E_iV_0$ have dimension 1 for $0 \leq i \leq D$ \cite[Lemma 
3.6]{terwSub1}. For a detailed description of $V_0$ see
\cite{curtin1}, \cite{egge1}. 

\bigskip \noindent
For the rest of this section assume $\Gamma$ is bipartite.  
There exists, up to isomorphism, a unique irreducible $T$-module with endpoint 1 
\cite[Corollary 7.7]{curtin1}.  We call this module $V_{1}$.  
The module $V_1$ is thin; in fact each of
$E^*_iV_1$,
$E_iV_1$ has dimension 1 for $1 \leq i \leq D-1$ and 
$E^{*}_{D}V_{1}=0$, $E_{0}V_{1}=0$, $E_{D}V_{1}=0$. 
For a detailed description of $V_1$ see
\cite{curtin1}. In this paper we are concerned with the 
thin irreducible $T$-modules with endpoint 2.

\bigskip \noindent
In order to describe the thin irreducible $T$-modules with endpoint 2 we define 
some parameters.
Let $\Gamma_{2}^{2} = \Gamma_{2}^{2}(x) $ denote the
graph with vertex set $\breve{X}$ and edge set $\breve{R}$, where
\begin{eqnarray}
\breve{X} &=& \{y \in X\;|\;\partial(x,y)=2\},\nonumber\\
\breve{R} &=& \{yz\;|\;y,z \in \breve{X},\, \partial(y,z)=2 \},\nonumber
\end{eqnarray}
and where $\partial$ is the path-length distance function for $\Gamma$.
The graph $\Gamma_{2}^{2}$ has 
exactly $k_{2}$ vertices, where $k_{2}$ is the second valency of
$\Gamma.$ Also, $\Gamma_{2}^{2}$ is regular with valency $p^{2}_{22}$.
We let
$\eta_1, \eta_2 , \ldots , \eta_{k_{2}}$
denote the eigenvalues of the adjacency matrix of $\Gamma_{2}^{2}$. 
By \cite[Theorem 11.7]{curtin2}, these eigenvalues may be ordered such that $\eta_{1} = p_{22}^{2}$ 
and $\eta_{i}= b_{3}-1$ $(2 \leq i \leq k)$.   

\bigskip \noindent
Abbreviate $d = \lfloor D/2 \rfloor$.  It is shown in \cite[Theorem 11.4]{maclean5} that 
${\tilde \theta}_1 \leq 
\eta_i \leq {\tilde \theta}_d $ for $k+1 \leq i \leq k_{2}$, where
${\tilde \theta}_1=-1-b_2 b_{3}(\theta_1^{2}-b_{2})^{-1}$ and 
${\tilde \theta}_d=-1-b_2 b_{3}(\theta_d^{2}-b_{2})^{-1}$. 
We remark 
$\theta_{1}^{2} > b_{2} > \theta_{d}^{2}$ by \cite[Lemma 2.6]{jkt},
so ${\tilde \theta}_1 < -1$ and ${\tilde \theta}_d\geq 0$.

\bigskip \noindent
Let $W$ denote a thin irreducible $T$-module with endpoint 2.
Observe $E^*_2W$ is a $1$-dimensional eigenspace for $E^*_2A_{2}E^*_2$;
let $\eta $ denote the corresponding eigenvalue.
It turns out $\eta $
is among $\eta_{k+1}, \eta_{k+2},\ldots, \eta_{k_{2}}$ so
${\tilde \theta}_1 \leq \eta \leq {\tilde \theta}_d$.
We call $\eta $ the {\it local eigenvalue } of $W$.
To describe the structure of $W$ we distinguish four cases: (i) $\eta 
={\tilde \theta}_1 $; (ii) $D$ is odd
and $\eta = {\tilde \theta}_{d}$; (iii) $D$ is even 
and $\eta = {\tilde \theta}_{d}$; (iv) $ {\tilde \theta}_1 < \eta < {\tilde 
\theta}_d$.  
In \cite{maclean5} we investigated cases (i), (ii).
In the present paper we investigate cases (iii), (iv).

\bigskip \noindent
Concerning cases (i), (ii) our results from \cite{maclean5} are summarized as follows.
Choose $n \in \{1,d\}$ if $D$ is odd, and let $n=1$ if $D$ is even.  
Define $\eta = {\tilde \theta}_{n}$.  
Let $W$ denote a thin irreducible
$T$-module with endpoint 2 and local eigenvalue $\eta$. 
Then $W$ has dimension $D-3$.   
Let $v$ denote a nonzero vector in $E^*_2W$.
We showed $W$ has a basis
$E_iv$ $(1 \leq i \leq D-1,\;i\not=n,\; i \not= D-n)$.
We showed this basis is orthogonal (with respect to the Hermitian dot
product) and we computed the square-norm of each basis vector.
We showed $W$ has a basis
$E^*_{i+2}A_iv$ $(0 \leq i \leq D-4)$.
We found the matrix representing $A$ with respect to this basis.
We showed this basis is orthogonal and we computed the square-norm of 
each basis vector. We found the transition matrix relating our two 
bases for $W$.
We showed the following scalars are equal: (i) The multiplicity with 
which $W$ appears
in the standard module $\mathbb{C}^{X}$; (ii) The number of times $\eta$ 
appears among
$\eta_{k+1}, \eta_{k+2}, \ldots, \eta_{k_{2}}$.

\bigskip \noindent
Concerning case (iii) above, in the present paper we obtain the following results.
Assume $D$ is even, and let $W$ denote a thin irreducible $T$-module with endpoint 2 and local eigenvalue $ {\tilde \theta}_d$.
We show 
 the dimension of $W$ is $D-2$.
Let $v$ denote a nonzero vector in $E^*_2W$.
We show $W$ has a basis $E_iv $ $(1\leq i \leq D-1, \; i \not= d)$.
We show this basis is orthogonal and 
we compute the square-norm of each basis vector. 
We show $W$ has a basis   
$E^*_{i+2}A_iv$ $(0 \leq i \leq D-3)$.
We  find 
the matrix representing $A$ with respect to this basis.
We show this basis is orthogonal 
and we compute the square-norm of each basis vector. 
We find the transition matrix relating our two bases for $W$.

\bigskip \noindent
Concerning case (iv) above, in the present paper we obtain the following results.
Let $W$ denote a thin irreducible $T$-module with endpoint 2 and local eigenvalue $\eta$ $({\tilde \theta}_1 < \eta < {\tilde \theta}_d)$.
We show 
 the dimension of $W$ is $D-1$.
Let $v$ denote a nonzero vector in $E^*_2W$.
We show $W$ has a basis $E_iv $ $(1\leq i \leq D-1)$.
We show this basis is orthogonal and 
we compute the square-norm of each basis vector. 
We show $W$ has a basis   
$E^*_{i+2}A_iv$ $(0 \leq i \leq D-2)$.
We  find 
the matrix representing $A$ with respect to this basis.
We show this basis is orthogonal 
and we compute the square-norm of each basis vector. 
We find the transition matrix relating our two bases for $W$.

\medskip
\noindent 
For all $\eta \in \R$ let $\mu_\eta $ denote the multiplicity
with which $W$ appears in ${\C}^X$, where
$W$ is a thin irreducible $T$-module with endpoint 2 and local
eigenvalue $\eta$. If no such $W$ exists we interpret 
 $\mu_\eta=0$. 
We show 
 $\mu_\eta$  is at most  the number of times $\eta $ appears
 among $\eta_{k+1}, \eta_{k+2},\ldots, \eta_{k_{2}}$.  
Concerning the case of equality,
we show the following are equivalent: (i) For all $\eta \in \R$,  
 $\mu_\eta$  is equal to the number of times $\eta $ appears
 among $\eta_{k+1}, \eta_{k+2},\ldots, \eta_{k_{2}}$;  
(ii) Every irreducible $T$-module with endpoint 2 is thin.


\section{Preliminaries concerning distance-regular graphs}
In this section 
we review some definitions and basic concepts concerning distance-regular
graphs.
For more background information we refer the reader to 
\cite{bannai}, \cite{bcn}, \cite{godsil} or \cite{terwSub1}.

\medskip
\noindent
Let $X$ denote a nonempty  finite  set.
Let
 $\hbox{Mat}_X(\C)$ 
denote the $\C$-algebra
consisting of all matrices whose rows and columns are indexed by $X$
and whose entries are in $\;\C  $. Let
$V=\C^X$ denote the vector space over $\C$
consisting of column vectors whose 
coordinates are indexed by $X$ and whose entries are
in $\C$.
We observe
$\hbox{Mat}_X(\C)$ 
acts on $V$ by left multiplication.
We endow $V$ with the Hermitian inner product $\langle \, , \, \rangle$ 
which satisfies
$\langle u,v \rangle = u^t\overline{v}$ for all $u,v \in V$,
where $t$ denotes transpose and $-$ denotes complex conjugation.
We abbreviate $\|u\|^2 = \langle u,u \rangle$ for all 
$u \in V.$ 
For all $y \in X,$ let $\hat{y}$ denote the element
of $V$ with a 1 in the $y$ coordinate and 0 in all other coordinates.
We observe $\{\hat{y}\;|\;y \in X\}$ is an orthonormal basis for $V.$
The following formula will be useful. For all  
$B \in \hbox{Mat}_X(\C)$ and for all $u,v \in V$, 
\begin{equation}\label{ADJ}
\langle Bu,v \rangle =  
\langle u, {\overline B}^tv \rangle. 
\end{equation}

\medskip
\noindent
Let $\Gamma = (X,R)$ denote a finite, undirected, connected graph,
without loops or multiple edges, with vertex set $X$ and 
edge set
$R$.   
Let $\partial $ denote the
path-length distance function for $\Gamma $,  and set
$D = {\rm max}\{\partial(x,y) \;|\; x,y \in X\}. $
We refer to $D$ as the {\it diameter} of $\Gamma $.
Let $\lfloor D/2 
\rfloor$ denote the greatest integer at most $D/2$. 
Vertices $x,y \in X$ are called 
{\it adjacent} whenever $xy$ is an edge.
For an integer $k \geq 0$, we say $\Gamma $ is 
 {\it regular} with {\it valency k} whenever each vertex of 
$\Gamma$ is adjacent to exactly $k$ distinct vertices of $\Gamma $. 
 We say $\Gamma$ is {\it distance-regular}
whenever for all integers $h,i,j\;(0 \le h,i,j \le D)$ 
and for all
vertices $x,y \in X$ with $\partial(x,y)=h,$ the number
\begin{equation}\label{PHIJ}
p_{ij}^h = |\{z \in X \; |\; \partial(x,z)=i, \partial(z,y)=j \}|
\end{equation}
is independent of $x$ and $y.$ The $p_{ij}^h$ are called
the {\it intersection numbers} of $\Gamma.$ We abbreviate
$c_i= p_{1i-1}^i \;(1 \le i \le D),\;a_i = p_{1i}^i \; (0 \le i \le D)$,
and 
$b_i= p_{1i+1}^i \; (0 \le i \le D-1)$.  For notational
convenience, we define $c_0=0$ and $b_D=0.$ We note 
$a_0=0$ and $c_1=1$. 

\medskip
\noindent
For the rest of this paper we assume  $\Gamma$  
is  distance-regular  with diameter $D\geq 3$. 

\medskip
\noindent
By (\ref{PHIJ}) and the triangle inequality,
\begin{equation}\label{REL1}
p_{1j}^h = 0 \qquad {\rm if} \qquad |h-j| > 1 \qquad \qquad (0 \le h,j \le D).
\end{equation}
Observe $\Gamma$ is regular with valency $k=b_0,$ and that
$c_i+a_i+b_i = k$ for $0 \le i \le D$.
Moreover $b_i > 0\;(0 \le i \le D-1)$ and $c_i>0 \;(1 \le i \le D).$
For $0 \le i \le D$ we abbreviate $k_i=p_{ii}^0,$ and observe
\begin{equation}\label{DEFKI}
k_i = |\{z \in X \;|\; \partial(x,z)=i\}|,
\end{equation}
where $x$ is any vertex in $X$. Apparently $k_0=1$ and $k_1 =k.$
By \cite[p.195]{bannai} we have
\begin{equation}\label{KI}
k_i = \frac{b_0b_1\cdots b_{i-1}}{c_1c_2\cdots c_i} \qquad \qquad (0 \le i \le D).
\end{equation}

\medskip
\noindent
We recall the Bose-Mesner algebra of $\Gamma.$ 
For 
$0 \le i \le D$ let $A_i$ denote the matrix in $\hbox{Mat}_X(\C)$ with
$xy$ entry
$$
{(A_i)_{xy} = \cases{1, & if $\partial(x,y)=i$\cr
0, & if $\partial(x,y) \ne i$\cr}} \qquad (x,y \in X).
$$
We call $A_i$ the $i^{\hbox{th}}$ {\it distance matrix} of $\Gamma.$
For convenience we define $A_i=0$ for $i < 0$ and $i > D.$
We abbreviate $A=A_1$ and call this the {\it adjacency
matrix} of $\Gamma.$ We observe
(ai) $A_0 = I$; 
(aii) $\sum_{i=0}^D A_i = J$;
(aiii) $\overline{A}_i = A_i $ $(0 \le i \le D)$;
(aiv) $A_i^t = A_i $ $(0 \le i \le D)$;
(av) $A_iA_j = \sum_{h=0}^D p_{ij}^h A_h$ $( 0 \le i,j \le D)$,
where $I$ denotes the identity matrix 
and $J$ denotes the all 1's matrix.
Let $M$ denote the subalgebra of $\hbox{Mat}_X(\C)$ generated by $A.$ 
Using
(ai), (av)
 one can readily show  $A_0,A_1,\ldots,A_D$
form a basis for $M.$ 
We refer to $M$ as the {\it Bose-Mesner algebra} of $\Gamma$.
By \cite[p.45]{bcn} $M$ has a second basis 
$E_0,E_1,\ldots,E_D$ such that
(ei) $E_0 = |X|^{-1}J$;
(eii) $\sum_{i=0}^D E_i = I$;
(eiii) $\overline{E}_i = E_i$ $(0 \le i \le D)$;
(eiv) $E_i^t = E_i$ $(0 \le i \le D)$;
(ev) $E_iE_j = \delta_{ij}E_i$ $(0 \le i,j \le D)$.
We refer to $E_0, E_1, \ldots, E_D $ as the {\it primitive idempotents}
of $\Gamma$.  
We  call $E_0$ the {\it trivial idempotent} of $\Gamma.$

\medskip
\noindent 
We  recall the eigenvalues
of  $\Gamma $.
Since $E_0,E_1,\ldots,E_D$ form a basis for  
$M,$ there exist complex scalars $\theta_0,\theta_1,
\ldots,\theta_D$ such that
$A = \sum_{i=0}^D \theta_iE_i$.
Combining this with 
(ev)
we find 
$AE_i = E_iA =  \theta_iE_i$ for $0 \leq i \leq D$.
Using (aiii) and
(eiii) we find $\theta_0,\theta_1,\ldots,\theta_D$ are in $\R.$ Observe
$\theta_0,\theta_1,\ldots,\theta_D$ are distinct 
since $A$ generates $M.$ By 
\cite[Proposition 3.1]{biggs} we have 
$\theta_0=k$
and $-k \le \theta_i \le k$ for $0 \le i \le D.$ Throughout this paper
we assume $E_0, E_1,\ldots, E_D$ are indexed so that 
$\theta_0 > \theta_1 >\cdots >\theta_D.$ 
We refer to $\theta_i$ as the {\it eigenvalue} of $\Gamma$ associated
with $E_i.$ We call $\theta_0$ the {\it trivial eigenvalue} of $\Gamma.$
For $0 \leq i \leq D$ let $m_i$ denote the rank of $E_i$. 
We refer to $m_i$ as the {\it multiplicity} of $E_i$ (or $\;\theta_i$).
From (ei) we find 
$m_0=1$. 
Using (eii)--(ev)
we find 
\begin{equation}\label{EIVDECOM}
V = E_0V+E_1V+ \cdots +E_DV \qquad \qquad {\rm (orthogonal\ direct\ sum}).
\end{equation}
For $0 \le i \le D$ the space $E_iV$ is the eigenspace of $A$ associated 
with $\theta_i$. We observe the dimension of $E_iV$ is $m_i$.
We now record a fact about the eigenvalues 
 $\theta_1$ and $\theta_D$.

\begin{lemma}\cite[Lemma 2.6]{jkt}\label{THETAJ}
Let $\Gamma $ denote a distance-regular graph with diameter $D \ge 3$ 
and eigenvalues $k=\theta_0 > \theta_1 > \cdots >\theta_D.$ Then
(i) 
 $-1 < \theta_1 < k$; (ii)
$a_1-k \le \theta_D < -1$.
\end{lemma}
Later in this paper we will discuss polynomials in one or two
variables. We will use the following notation.
Let  $\lambda $ 
denote an indeterminate. Let
 $\R \lbrack \lambda \rbrack $ denote the 
$\R$-algebra consisting of all polynomials in $\lambda $ that
have coefficients in $\R$.
Let $\mu $ denote an indeterminate 
which commutes with $\lambda $.  Let  
 $\R \lbrack \lambda, \mu \rbrack $ denote the 
$\R$-algebra consisting of all polynomials in $\lambda $ and $\mu$ that
have coefficients in $\R$. 


\section{Bipartite distance-regular graphs}

\bigskip

        We now consider the case in which $\Gamma$ is bipartite.  
        We say $\Gamma$ is {\it bipartite}
        whenever the vertex set $X$ can be partitioned into two subsets, neither
        of which contains an edge. In the next few lemmas, we recall some routine
        facts concerning the case in which $\Gamma$ is bipartite.  
        To avoid trivialities, we will 
        generally assume $D \geq 4$.   
        
\begin{lemma} \label{bipequiv}  \cite[Propositions 3.2.3, 4.2.2]{bcn} \quad
        Let $\Gamma$ denote a distance-regular graph with diameter $D \geq 4$,
        valency $k$, 
        and eigenvalues
        $\theta_0 > \theta_1 > \cdots > \theta_D$.  The following are
        equivalent:
\begin{description}
\item[{\rm (i)}] 
       $\Gamma$ is bipartite.
\item[{\rm (ii)}]
       $p^{h}_{ij} = 0 \mbox{ if } h+i+j \mbox{ is odd} \qquad (0 \le 
       h,i,j \le D)$.
\item[{\rm (iii)}]
        $a_i = 0 \qquad (0 \le i \le D)$.
\item[{\rm (iv)}]
        $c_i + b_i =k  \qquad (0 \le i \le D)$.
\item[{\rm (v)}]
        $\theta_{D-i} = -\theta_i \qquad (0 \le i \le D).$
\end{description} 
\end{lemma}

\begin{lemma}  \label{thetad}
Let $\Gamma$ denote a bipartite distance-regular graph with diameter $D \geq 4$
        and eigenvalues
        $k=\theta_0 > \theta_1 > \cdots > \theta_D$.
\begin{description}
\item[{\rm (i)}] Assume $D$ is even and let $d = D/2$.  Then $\theta_{d}=0$.
\item[{\rm (ii)}] Assume $D$ is odd and let $d =(D-1)/2$.  Then 
$\theta_{d} > 0$ and $\theta_{d+1}=-\theta_{d}$.
\end{description}
\end{lemma}        
{\it Proof.}  Immediate from Lemma \ref{bipequiv}(v).  \hfill $\Box$\\

\bigskip

\begin{lemma}  \cite[Lemma 3.4]{maclean5}   \label{J-ED}
   Let $\Gamma =(X,R)$ denote a bipartite distance-regular graph with diameter 
   $D \geq 4$
        and eigenvalues
        $\theta_0 > \theta_1 > \cdots > \theta_D$.   Then $E_{D} = 
        |X|^{-1}J'$, where 
\begin{equation}  \label{J-def}
     J' = \sum_{i=0}^{D} (-1)^{i}A_{i}.
\end{equation}
\end{lemma}

\begin{lemma}  \label{bipeigs}
Let $\Gamma$ denote a bipartite distance-regular graph with diameter 
$D \geq 4$ and eigenvalues $\theta_{0} > \theta_{1} > \cdots > 
\theta_{D}$.  Then $\theta_{1}^{2} > 
b_{2} > \theta_{d}^{2},$ where $d = \lfloor D/2 \rfloor$.
\end{lemma}
{\it Proof.}  Apply Lemma \ref{THETAJ} to the halved graph of 
$\Gamma$, and use \cite[Proposition 4.2.3]{bcn}.
\hfill $\Box$\\



\section{Two families of polynomials}

\bigskip
\noindent
Let $\Gamma=(X,R)$ denote a bipartite distance-regular graph with diameter $D 
\ge 4.$ In this section we recall two types of polynomials associated 
with
$\Gamma $. To motivate things, we recall by 
(av) and the triangle inequality that \begin{equation}\label{AAI1}
AA_i = b_{i-1}A_{i-1}  + c_{i+1}A_{i+1} \qquad \qquad
(0 \le i \le D),
\end{equation}
where $b_{-1}=0$ and $c_{D+1}=0$.
Let $f_0,f_1,\ldots,$$f_D$ denote the polynomials in $\R[\lambda]$ 
satisfying $f_0=1$ and
\begin{equation}\label{FIPOLY2}
\lambda f_i = b_{i-1}f_{i-1}  + c_{i+1}f_{i+1} \qquad (0 \le 
i \le D-1),
\end{equation}
where $f_{-1}=0$.
For $0 \leq i \leq D$ the polynomial $f_i$ 
has degree $i$, and
the coefficient of
$\lambda^i$ is $(c_1c_2\cdots c_i)^{-1}.$
Comparing (\ref{AAI1}) and (\ref{FIPOLY2}) we find
$f_i(A)=A_i$. By \cite[p. 63]{bannai}
the polynomials $f_0,f_1,\ldots, f_D$ satisfy the orthogonality 
relation
\begin{displaymath}
\sum_{h=0}^D f_i(\theta_h)f_j(\theta_h)m_h = \delta_{ij} |X|k_i
\qquad \qquad (0 \le i,j \le D).
\end{displaymath}

\bigskip
\noindent We now recall some polynomials related to the $f_i$. 
Let $p_0, p_1, \ldots, p_D$ denote the polynomials in ${\R}[\lambda]$ 
satisfying 
\begin{equation}
{p_{i} = \cases{f_{0}+f_{2}+f_{4}+\cdots +f_{i}, & if $i$ is even \cr
f_{1}+f_{3}+f_{5}+\cdots +f_{i}, & if $i$ is odd \cr}} 
 \qquad \qquad (0 \leq i \leq D).
\label{PIPOLY}
\end{equation}

\bigskip
\noindent
Observe $p_{0} =1$.  
For $0 \leq i \leq D$
the polynomial $p_i$ has degree $i$, and the
coefficient of $\lambda^i$ is
$(c_1c_2\cdots c_i)^{-1}$.  Recalling $f_{j}(A)=A_{j} \; (0 \le j \le 
D)$, we observe

\begin{equation}  \label{PDAeqs}
     p_{D}(A) + p_{D-1}(A) = J, \qquad \qquad p_{D}(A) - 
     p_{D-1}(A) = (-1)^{D}J',
\end{equation}
where $J'$ is from (\ref{J-def}).  By \cite[Theorem 4.2]{maclean5}, we 
have 
\begin{equation}  \label{PIRECUR1}
\lambda p_i = c_{i+1}p_{i+1}+b_{i+1}p_{i-1}
\qquad (0 \leq i \leq D-1),
\end{equation}
where $p_{-1}=0$.  
We record a fact for later use.

\begin{lemma}\label{PIPJ}  \cite[Lemma 4.3]{maclean5}
Let $\Gamma=(X,R)$ denote a bipartite distance-regular graph with diameter $D 
\ge 4$ and eigenvalues $k=\theta_0 >\theta_1 >\cdots >\theta_D.$ Let 
the polynomials $p_0,p_1,\ldots, p_D$ be as in (\ref{PIPOLY}).
Then $p_{D-1}(\theta_h)=0$ and $p_{D}(\theta_h)=0$  for $1 \leq h \leq D-1$.
Moreover,
\begin{equation}
\label{PIPJ1}
\sum_{h=0}^D p_i(\theta_h)p_j(\theta_h)(k^{2}-\theta_h^{2})m_h = \delta_{ij} 
|X|k_ib_ib_{i+1}
\qquad \qquad (0 \le i,j \le D-2).
\end{equation}
\end{lemma}


\section{The polynomials $\Psi_{i}$}

\bigskip
\noindent  Let $\Gamma$ denote a bipartite distance-regular graph with 
diameter $D \geq 4$.  In the previous section we used $\Gamma$ to 
define two families of polynomials in one variable.  We called these 
polynomials the $f_{i}$ and the $p_{i}$.  Later in this paper we will 
use $\Gamma$ to define a third family of polynomials in one 
variable.  We will call these polynomials the $g_{i}$.  To define and 
study the $g_{i}$ it is convenient to first consider some 
polynomials $\Psi_{i}$ in two variables.

\begin{definition} \label{PSIPOLY}  \rm
Let $\Gamma$ denote a bipartite distance-regular graph with diameter $D 
\ge 4$.  For $0 \leq i \leq D-2$ let $\Psi_{i}$ denote the polynomial 
in $\R[\lambda, \mu]$ given by
\begin{equation}
\Psi_i = \sum_{{h=0}\atop{ i-h { \mbox{ \tiny even}}}}^{i} 
p_{h}(\lambda) p_{h}(\mu)
\frac{k_{i}b_{i}b_{i+1}}
{k_{h}b_{h}b_{h+1}},
\label{PSIPOLY1}
\end{equation}
where the polynomials $p_0, p_{1}, \ldots, p_{D-2}$ are from
(\ref{PIPOLY}).  We observe $\Psi_{0}=1$ and $\Psi_{1} = \lambda 
\mu$.  
\end{definition}

\begin{lemma} 
Let $\Gamma$ denote a bipartite distance-regular graph with diameter 
$D \geq 4$.  Let the polynomials $p_{i}, \Psi_{i}$ be as in 
(\ref{PIPOLY}), (\ref{PSIPOLY1}), respectively.  Then
\begin{displaymath}  
   p_{i}(\lambda) p_{i}(\mu) = \Psi_{i} 
   -\frac{b_{i}b_{i+1}}{c_{i}c_{i-1}}\Psi_{i-2} \qquad  \qquad (2 
   \leq i \leq D-2).
\end{displaymath}
\end{lemma}
{\it Proof.}  Use Definition \ref{PSIPOLY} and (\ref{KI}).  \hfill $\Box$\\

\bigskip \noindent
The following equation is a variation of the Christoffel-Darboux 
formula.

\begin{lemma}  \label{CDlemma}  \cite[Lemma 5.3]{maclean5}
Let $\Gamma$ denote a bipartite distance-regular graph with diameter 
$D \geq 4$.  Let the polynomials $p_{i}, \Psi_{i}$ be as in 
(\ref{PIPOLY}), (\ref{PSIPOLY1}) respectively.  Then for $1 \leq i 
\leq D-1$,
\begin{displaymath} 
  p_{i+1}(\lambda) p_{i-1}(\mu) - p_{i-1}(\lambda)p_{i+1}(\mu) = 
  c_{i}^{-1} c_{i+1}^{-1}(\lambda^{2}-\mu^{2})\Psi_{i-1}.
\end{displaymath}
\end{lemma}

\begin{lemma}  \label{PSIIPSIJ}  \cite[Lemma 5.4]{maclean5}
Let $\Gamma =(X,R)$ denote a bipartite distance-regular graph with diameter $D 
\ge 4$ and eigenvalues $k=\theta_0 >\theta_1 >\cdots >\theta_D.$  Let 
the polynomials $p_{i}, \; \Psi_{i}$ be as in (\ref{PIPOLY}), 
(\ref{PSIPOLY1}) respectively.  Then for $0 \leq i,j \leq D-2$,
\begin{displaymath} 
\sum_{h=0}^{D}\Psi_{i}(\theta_{h},\mu)\Psi_{j}(\theta_{h},\mu)(k^{2}-\theta_{h}^{2})
(\mu^{2}-\theta_{h}^{2})m_{h} 
  \; = \; \delta_{ij}|X|p_{i}(\mu)p_{i+2}(\mu)k_{i}b_{i}b_{i+1}c_{i+1}c_{i+2}.
 \end{displaymath}
 (We recall $m_{h}$ denotes the multiplicity of $\theta_{h}$ for $0 
 \leq h \leq D.$)
\end{lemma}

\begin{lemma}  \label{PITHETA}
Let $\Gamma$ denote a bipartite distance-regular graph with diameter $D 
\ge 4$ and eigenvalues $k=\theta_0 >\theta_1 >\cdots >\theta_D.$  Let 
the polynomials $p_{i}$ be as in (\ref{PIPOLY}).  Then 
the following (i), (ii) hold 
for all $\theta \in \R$:
\begin{description}
\item[{\rm (i)}]  Suppose $\theta=\theta_{1}$.  Then $p_{i}(\theta) > 0 $
for $0 \le i \le D-2$, and $p_{D-1}(\theta)=0, \, p_{D}(\theta)=0$.
\item[{\rm (ii)}]  Suppose $\theta > \theta_{1}$.  Then $p_{i}(\theta) > 0 $
for $0 \le i \le D$. 
\end{description}
\end{lemma}
{\it Proof.}   Observe $p_{D-1}(\theta_{1})=0, \, 
p_{D}(\theta_{1})=0$ by Lemma \ref{PIPJ}.  For notational convenience 
set $e=0$ if $\theta > \theta_{1}$ and $e=1$ if $\theta = 
\theta_{1}$.  Suppose there exists an integer $i$ $(0 \le i \le D-2e)$ 
such that $p_{i}(\theta) \le 0$.  Let us pick the minimal such $i$.  
Observe $i \ge 2$ since $p_{0}(\theta) =1$, $p_{1}(\theta) = 
\theta$.  Apparently $p_{i-2}(\theta) > 0$.  We claim there exists an 
integer $h$ $(1+e \le h \le D-1-e)$ such that $\Psi_{i-2}(\theta_{h}, 
\theta) \not= 0$.  To see this, observe by Definition \ref{PSIPOLY} 
that $\Psi_{i-2}(\lambda, \theta)$ is a polynomial in $\lambda$ with 
degree $i-2$.  In this polynomial the coefficient of $\lambda^{i-2}$ 
is $p_{i-2}(\theta)(c_{1}c_{2}\cdots c_{i-2})^{-1}$.  Apparently this 
polynomial is not identically 0 so there exist at most $i-2$ integers 
$h$ $(1+e \le h \le D-1-e)$ such that $\Psi_{i-2}(\theta_{h}, 
\theta)=0$.  By this and since $i \le D-2e$, there exists at least one 
integer $h$ $(1+e \le h \le D-1-e)$ such that $\Psi_{i-2}(\theta_{h}, 
\theta) \not= 0$.  We have now proved our claim.  We may now argue
\begin{eqnarray*}
 0 &<& \sum_{h=1+e}^{D-1-e} 
 \Psi_{i-2}^{2}(\theta_{h},\theta)(k^{2}-\theta_{h}^{2}) 
 (\theta^{2}-\theta_{h}^{2})m_{h} \\
 &=& \sum_{h=0}^{D} 
 \Psi_{i-2}^{2}(\theta_{h},\theta)(k^{2}-\theta_{h}^{2}) 
 (\theta^{2}-\theta_{h}^{2})m_{h} \quad \mbox{ (by the definition of $e$)} \\
 &=& |X|p_{i-2}(\theta)p_{i}(\theta)k_{i-2}b_{i-2}b_{i-1}c_{i-1}c_{i} 
 \quad
 \mbox{ (by Lemma \ref{PSIIPSIJ}) }\\ & \le& 0.
 \end{eqnarray*}
 We now have a contradiction and the result follows. 
 \hfill $\Box$\\

\begin{lemma}  \label{PITHETAD}
Let $\Gamma$ denote a bipartite distance-regular graph with odd diameter $D 
\ge 4$ and eigenvalues $k=\theta_0 >\theta_1 >\cdots >\theta_D.$  
Let $d$ denote the integer satisfying $2d+1 =D$.  Let 
the polynomials $p_{i}$ be as in (\ref{PIPOLY}).  Then the following (i), (ii) hold 
for all $\theta \in \R$:
\begin{description}
\item[{\rm (i)}]  Suppose $\theta=\theta_{d}$.  Then $(-1)^{\lfloor 
\frac{i}{2}\rfloor}p_{i}(\theta) > 0 $
for $0 \le i \le D-2$, and $p_{D-1}(\theta)=0, \, p_{D}(\theta)=0$.
\item[{\rm (ii)}]  Suppose $0 < \theta < \theta_{d}$.  Then $(-1)^{\lfloor 
\frac{i}{2}\rfloor}p_{i}(\theta) > 0 $
for $0 \le i \le D$. 
\end{description}
\end{lemma}
{\it Proof.}  Observe $p_{D-1}(\theta_{d})=0, \, 
p_{D}(\theta_{d})=0$ by Lemma \ref{PIPJ}.  For notational convenience 
set $e=0$ if $0 < \theta < \theta_{d}$ and $e=1$ if $\theta = 
\theta_{d}$.  Also for notational convenience we define the set $S$ to be $\{1, 2, \ldots, D-1\}$ if $e=0$, and $\{1, 2, \ldots, d-1\} \cup \{d+2, d+3, \ldots, D-1\}$ if $e=1$.
Suppose there exists an integer 
$i$ $(0 \le i \le D-2e)$ 
such that $(-1)^{\lfloor \frac{i}{2}\rfloor}p_{i}(\theta) \le 0$.  
Let us pick the minimal such $i$.  
Observe $i \ge 2$ since $p_{0}(\theta) =1$, $p_{1}(\theta) = 
\theta$.  Apparently $(-1)^{\lfloor \frac{i-2}{2}\rfloor}p_{i-2}(\theta) > 
0$, so $p_{i-2}(\theta) p_{i}(\theta) \geq 0$.    
We claim there exists an 
integer $h \in S$ such that 
$\Psi_{i-2}(\theta_{h}, 
\theta) \not= 0$.  To see this, observe by Definition \ref{PSIPOLY} 
that $\Psi_{i-2}(\lambda, \theta)$ is a polynomial in $\lambda$ with 
degree $i-2$.  This polynomial is not identically zero, since the coefficient of 
 $\lambda^{i-2}$ 
is $p_{i-2}(\theta)(c_{1}c_{2}\cdots c_{i-2})^{-1}$ and since $p_{i-2}(\theta) \not= 0$ by construction.
Therefore there exist at most $i-2$ integers $h \in S$
such that $\Psi_{i-2}(\theta_{h}, 
\theta)=0$.  By this and since $i \le D-2e$, there exists at least one 
integer $h \in S$ 
such that $\Psi_{i-2}(\theta_{h}, 
\theta) \not= 0$.  We have now proved our claim.  We may now argue
\begin{eqnarray*}
 0 &>& \sum_{h \in S} 
 \Psi_{i-2}^{2}(\theta_{h},\theta)(k^{2}-\theta_{h}^{2}) 
 (\theta^{2}-\theta_{h}^{2})m_{h}  \\
 &=& \sum_{h=0}^{D} 
 \Psi_{i-2}^{2}(\theta_{h},\theta)(k^{2}-\theta_{h}^{2}) 
 (\theta^{2}-\theta_{h}^{2})m_{h} \quad \mbox{ (by the definitions of $S$ and $e$)}  \\
 &=& |X|p_{i-2}(\theta)p_{i}(\theta)k_{i-2}b_{i-2}b_{i-1}c_{i-1}c_{i} 
 \quad
 \mbox{ (by Lemma \ref{PSIIPSIJ}) }\\ & \ge& 0.
 \end{eqnarray*}
 We now have a contradiction and the result follows. 
  \hfill $\Box$\\
 

\section{A variation of the $p_{i}$ polynomials}

In Section 4 we defined some polynomials $p_i$.  In this section we define some closely related polynomials that we call the $P_i$.  We do so for a technical reason that will become apparent later in the paper.  We start with an observation.  Recall that a polynomial in 
$\R[\lambda]$ is {\it even} (resp. {\it odd}) whenever the coefficient of $\lambda^i$ is zero for all odd $i$ (resp. all even $i$).

\begin{lemma} \label{p-even-odd} 
Let $\Gamma$ denote a bipartite distance-regular graph 
with diameter $D \ge 4$.  Then for $0 \leq i \leq D$ the polynomial $p_i$ from 
(\ref{PIPOLY}) is even (resp. odd)
if $i$ is even (resp. odd).
\end{lemma}
{\it Proof.}  Routine using (\ref{PIRECUR1}) and induction.  \hfill $\Box$\\

\medskip

In view of Lemma \ref{p-even-odd} we can make the following definition.

\begin{definition} \label{Pidef}  \rm
Let $\Gamma$ denote a bipartite distance-regular graph 
with diameter $D \ge 4$.  For $0 \leq i \leq D$ let $P_i$ denote the polynomial in $\R[\lambda]$ such that
\begin{equation}  \label{Pidef2}
{p_{i}(\lambda) = \cases{P_{i}(\lambda^{2}), & if $i$ is even \cr
\lambda P_{i}(\lambda^{2}), & if $i$  is odd, \cr}} 
\end{equation}
where $p_i$ is from (\ref{PIPOLY}).
Observe the degree of $P_{i}$ is $i/2$ if $i$ is even
and  $(i-1)/2$ if $i$ is odd.  For notational convenience we define $P_{-1}=0$.
\end{definition}

\begin{lemma} \label{Precur}
Let $\Gamma$ denote a bipartite distance-regular graph with diameter $D 
\ge 4.$   Let the polynomials $P_{0}, P_{1}, \ldots, 
P_{D}$ be as in Definition \ref{Pidef}.  
Then the following (i), (ii) hold for $0 \leq i \leq D-1$:
\begin{description}
\item[{\rm (i)}]  Suppose $i$ is odd.  Then $\lambda P_{i} = 
c_{i+1} P_{i+1}  + b_{i+1}P_{i-1}$.
\item[{\rm (ii)}]  Suppose $i$ is even.  Then $P_{i} = 
c_{i+1} P_{i+1} + b_{i+1}P_{i-1}$.  
\end{description}
\end{lemma}
{\it Proof.}   Routine using (\ref{PIRECUR1}) and Definition 
\ref{Pidef}.   \hfill $\Box$\\

Referring to Lemma \ref{Precur}, in order to handle the cases of $i$ odd and $i$ even in a uniform fashion we introduce some notation.

\begin{definition} \label{s(i)}  \rm
For any integer $i$ we define 
$$
s(i) =   \cases{0, & if $i$ is even \cr
1, & if $i$ is odd. \cr}
$$
\end{definition}

\noindent
Lemma \ref{Precur} looks as follows in terms of $s(i)$.

\begin{corollary} \label{Precur2}
Let $\Gamma$ denote a bipartite distance-regular graph with diameter $D 
\ge 4,$ and let the polynomials $P_{0}, P_{1}, \ldots, 
P_{D}$ be as in Definition \ref{Pidef}.  
Then for $0 \leq i \leq D-1,$
\begin{equation}  \label{Precur3}
\lambda^{s(i)} P_{i} = 
c_{i+1} P_{i+1}  + b_{i+1}P_{i-1}.
\end{equation}
\end{corollary}

\begin{lemma} \label{Pnot0}
Let $\Gamma$ denote a bipartite distance-regular graph with diameter $D 
\ge 4$  and eigenvalues $k=\theta_0 >\theta_1 >\cdots >\theta_D.$   
Let the polynomials $P_{0}, P_{1}, \ldots, 
P_{D}$ be as in Definition \ref{Pidef}.  Then the following (i)--(iii) hold for all $\psi \in \mathbb{R}$:
\begin{description}
\item[{\rm (i)}] Assume $\psi > \theta_1^2$.  Then $P_i(\psi) > 0$ $\; (0 \le i \le D)$.
\item[{\rm (ii)}] Assume $D$ is odd and $\psi < \theta_d^2$, where $d=(D-1)/2$.  Then $(-1)^{\lfloor \frac{i}{2} 
\rfloor}P_{i}(\psi) >0$ $\; (0 \le i \le D)$. 
\item[{\rm (iii)}]  Assume $D$ is even and $\psi \leq 0$.  Then $(-1)^{\lfloor \frac{i}{2} 
\rfloor}P_{i}(\psi) >0$ $\; (0 \le i \le D-1)$.   Moreover $(-1)^{\lfloor \frac{D}{2} 
\rfloor}P_{D}(\psi) >0$ if $\psi < 0$ and $P_D(0)=0$.
\end{description}
\end{lemma}
{\it Proof.}  (i).  Since $\psi$ is 
positive, there exists a positive real number $\alpha$ such that 
$\alpha^{2}= \psi$.  By the construction $\alpha > \theta_{1}$.  
For $0 \leq i \leq D$ we have $p_i(\alpha) >0$ by Lemma 
\ref{PITHETA}(ii) so 
$P_{i}(\psi) > 0$ in view of Definition \ref{Pidef}. \\
(ii).  First assume $0 < \psi < \theta_{d}^{2}$.  Again $\psi$ is 
positive, so there exists a positive real number $\alpha$ such that 
$\alpha^{2}= \psi$.  By the construction $ 0 < \alpha < \theta_{d}$.  
For $0 \leq i \leq D$ we have $(-1)^{\lfloor \frac{i}{2} \rfloor} p_{i}(\alpha) > 0$
by Lemma \ref{PITHETAD}(ii) so $(-1)^{\lfloor \frac{i}{2} \rfloor}P_{i}(\psi) > 0$ in view of Definition \ref{Pidef}.

Now assume $\psi \leq 0$.  Suppose there exists an integer $i$ $(0 \leq i \leq D)$ such that 
$(-1)^{\lfloor \frac{i}{2} \rfloor}P_{i}(\psi) \leq 0$.  Let us pick the minimal such $i$.  Observe $i \geq 2$ since $P_0(\psi) =1$, $P_1(\psi)=1$.  Setting $\lambda = \psi$ and replacing $i$ by $i-1$ in (\ref{Precur3}) and then multiplying this equation by $(-1)^{\lfloor \frac{i}{2} \rfloor}$, we find
\begin{eqnarray}
  (-1)^{\lfloor \frac{i}{2} \rfloor} P_{i}(\psi) &=&
    (-1)^{\lfloor \frac{i}{2} \rfloor}  \psi^{s(i-1)} P_{i-1}(\psi) c_{i}^{-1}
    -(-1)^{\lfloor \frac{i}{2} \rfloor} 
  P_{i-2}(\psi) b_{i}c_{i}^{-1}   \nonumber \\
  &=& (-\psi)^{s(i-1)}  (-1)^{\lfloor \frac{i-1}{2} \rfloor} P_{i-1}(\psi) c_{i}^{-1} 
  - (-1)^{\lfloor \frac{i}{2} \rfloor}P_{i-2}(\psi) b_{i}c_{i}^{-1} \label{warwick} \\
  & > & 0, \nonumber
\end{eqnarray}
where the last inequality follows from the minimality of $i$ and $\psi \leq 0$.  We now have a contradiction and the result follows. \\
(iii).  Similar to (ii).  When $\psi=0$, however, observe that the right side of (\ref{warwick}) is 0 for $i=D$, and hence $P_D(0)=0$.   \hfill $\Box$\\ 

\begin{corollary} \label{Pnot0v2}
Let $\Gamma$ denote a bipartite distance-regular graph with diameter $D 
\ge 4$ and eigenvalues $k=\theta_0 >\theta_1 >\cdots >\theta_D.$   
Let the polynomials $P_{0}, P_{1}, \ldots, 
P_{D}$ be as in Definition \ref{Pidef}.  Let $\theta$ denote a real number in the following range:
For $D$ odd, we assume  $\theta > \theta_{1}^{2}$ or 
$ \theta < \theta_{d}^{2}$, where $d=(D-1)/2$. 
For $D$ even, we assume $\theta > \theta_1^2$ or $\theta \leq 0$.
Then $P_{i}(\theta) \not= 0$ for $0 \le i \le 
D-1$.   
\end{corollary}


\section{A third family of polynomials}

\bigskip
\noindent In this section we will use the following notation.

\begin{notation} \label{thetanote}  \rm
Let $\Gamma =(X,R)$ denote a bipartite distance-regular graph with diameter $D \ge 4$ and eigenvalues $k=\theta_0 >\theta_1 >\cdots >\theta_D.$  
Let $d = \lfloor D/2 \rfloor$.  Let the polynomials $p_i$ be as in (\ref{PIPOLY}), and let
the polynomials $P_{i}$ be as in Definition \ref{Pidef}.
Let $\theta$ denote a real number in the following range:
For $D$ odd, we assume  $\theta > \theta_{1}^{2}$ or 
$ \theta < \theta_{d}^{2}$.
For $D$ even, we assume $\theta > \theta_1^2$ or $\theta \leq 0$.
  We observe that in all cases
 $P_{i}(\theta) \not=0$ for $0 \le i \le 
D-1$ by Corollary \ref{Pnot0v2}.
\end{notation}

\noindent We now use $\Gamma$
to define a family of polynomials in one variable.  We call these 
polynomials
the $g_i$.

\begin{definition} \label{GIPOLY}  \rm
With reference to Notation \ref{thetanote},
for $0 \leq i \leq D-2$ we define the polynomial
$g_i \in \R[\lambda]$ by
\begin{equation}
g_i = \sum_{{h=0}\atop{ i-h { \mbox{ \tiny even}}}}^{i}
\frac{P_h(\theta)}{P_i(\theta)}\frac{k_i b_i b_{i+1}}{ k_hb_h b_{h+1}} p_h.
\label{GIPOLY1}
\end{equation}  
We emphasize $g_i$ depends on $\theta$ as well as
the intersection numbers of $\Gamma$.
\end{definition}

\begin{lemma}  
With reference to Notation \ref{thetanote} and Definition \ref{GIPOLY}, 
\begin{equation}  \label{pigi}
  p_{i} = g_{i} - 
  \frac{b_{i}b_{i+1}}{c_{i-1}c_{i}}\frac{P_{i-2}(\theta)}{P_{i}(\theta)}g_{i-2} 
   \qquad (2 \le i \le D-2).
\end{equation} 
\end{lemma}
{\it Proof.}  Routine using Definition \ref{GIPOLY} and (\ref{KI}).  \hfill $\Box$\\

\begin{lemma}  \label{P8}
With reference to Notation \ref{thetanote}  and Definition \ref{GIPOLY}, 
the following 
 (i), (ii) hold for $0 \le i 
\le D-2$:
\begin{description}
\item[{\rm (i)}]  The polynomial $g_{i}$ has degree exactly $i$.  
\item[{\rm (ii)}]  The coefficient of $\lambda^{i}$ in $g_{i}$ is 
$(c_{1}c_{2}\cdots c_{i})^{-1}$.  
\end{description}
\end{lemma}
{\it Proof.}  Routine.  \hfill $\Box$\\

\bigskip \noindent
We now present a three-term recurrence satisfied by the polynomials $g_{i}$.

\begin{theorem}  \label{P9}  
With reference to Notation \ref{thetanote}  and Definition \ref{GIPOLY}, 
 $g_{0}=1$ and
\begin{equation}  \label{girecur}
  \lambda g_{i} = c_{i+1}g_{i+1} + \omega_{i} g_{i-1} 
\end{equation}
   for $0 \le i \le D-2$, where $g_{-1}=0, \,  \omega_{0}=0, \, 
   g_{D-1} = p_{D-1},$
and
\begin{equation}  \label{omega}
  \omega_{i} = \frac{b_{i+1}c_{i+2}}{c_{i}} \frac{P_{i-1}(\theta) 
  P_{i+2}(\theta)}{P_{i}(\theta) P_{i+1}(\theta)} \qquad \qquad (1 
  \le i \le D-2).
\end{equation}
\end{theorem}
{\it Proof.} We find $g_{0}=1$ by Definition \ref{GIPOLY}.  We now prove (\ref{girecur})
by induction on $i$.  Line (\ref{girecur}) holds for $i=0,1$ using Definition \ref{GIPOLY}, (\ref{PIRECUR1}), and Definition \ref{Pidef}.  Next 
assume $i \geq 2$
and by induction that  
\begin{equation}  \label{girecurA}
  \lambda g_{i-2} =  c_{i-1} g_{i-1} + \omega_{i-2} g_{i-3} .
\end{equation}
Consider the right-hand side of (\ref{girecur}).  
In this expression eliminate $g_{i+1}$ using (\ref{pigi}) if $i < D-2$ and $g_{D-1}=p_{D-1}$ if 
$i=D-2$.  Also eliminate $\omega_i$ using (\ref{omega}) and simplify the result using 
(\ref{Precur3}) to get
\begin{equation}   \label{girecurD}
  c_{i+1}g_{i+1} + \omega_{i} g_{i-1}
   = c_{i+1} p_{i+1} + \frac{b_{i+1} \theta^{s(i+1)} }{c_i}\frac{    P_{i-1}(\theta)}{
   P_{i}(\theta)} g_{i-1}.
\end{equation}
Now consider the left-hand side of (\ref{girecur}).  Replacing $g_{i}$ in 
this expression  using (\ref{pigi}), and eliminating $\lambda p_{i}$, 
$\lambda g_{i-2}$ in the result using (\ref{PIRECUR1}), 
(\ref{girecurA}), respectively, we find
\begin{equation}  \label{girecurB}
   \lambda g_{i} =  c_{i+1}p_{i+1} + b_{i+1} p_{i-1} + 
   \frac{b_{i}b_{i+1}}{c_{i-1}c_{i}} 
   \frac{P_{i-2}(\theta)}{P_{i}(\theta)} (c_{i-1} g_{i-1} + 
   \omega_{i-2} g_{i-3}). 
\end{equation}
If $i >2$, in (\ref{girecurB}) we eliminate $\omega_{i-2}$ using (\ref{omega}) and then eliminate 
$b_{i-1}b_{i}P_{i-3}(\theta)(c_{i-2}c_{i-1}P_{i-1}(\theta))^{-1} 
  g_{i-3}$ in the resulting equation using (\ref{pigi}).  If $i=2$, in (\ref{girecurB}) we note $\omega_0 =0$ and $p_1 =g_1$ in view of Definition \ref{GIPOLY}.  In either case we find
\begin{equation}  \label{girecurC}
 \lambda g_{i} = c_{i+1}p_{i+1} + b_{i+1} 
 \, \frac{c_{i}P_{i}(\theta) 
 +b_{i}P_{i-2}(\theta)}{c_{i}P_{i}(\theta)} \, g_{i-1}.
\end{equation}
Observe the right-hand sides of (\ref{girecurD}), (\ref{girecurC}) are 
equal
in view of (\ref{Precur3}) and Definition \ref{s(i)}, and thus the left-hand sides are equal.
We obtain (\ref{girecur}) as desired.
\hfill $\Box$\\

  \begin{lemma}
With reference to Notation \ref{thetanote}  and Definition \ref{GIPOLY}, 
 for $0 \le i \le D-2$ we have
\begin{equation}  \label{gipi}
  c_{i+1}^{-1}c_{i+2}^{-1}(\lambda^{2}-\theta) g_{i} = 
  p_{i+2}-\frac{P_{i+2}(\theta)}{P_{i}(\theta)}p_{i}.
\end{equation}
\end{lemma}
{\it Proof.}  We show (\ref{gipi}) by induction on $i$.  Line
(\ref{gipi}) holds for $i=0,1$ by Definition \ref{GIPOLY}, (\ref{PIRECUR1}), and Definition \ref{Pidef}.  Next assume $i \geq 2$ and by induction that
\begin{equation} \label{gipi0}
c_{i-1}^{-1}c_{i}^{-1} (\lambda^{2}-\theta)g_{i-2} = 
  p_{i}-\frac{P_{i}(\theta)}{P_{i-2}(\theta)}p_{i-2}.
\end{equation}
Repeatedly applying (\ref{PIRECUR1}), we find
\begin{equation}  \label{gipi1}
\lambda^{2} p_{i} = c_{i+1}c_{i+2}p_{i+2} + 
(c_{i+1}b_{i+2}+b_{i+1}c_{i}) p_{i} + b_{i}b_{i+1} p_{i-2}.
\end{equation}
Similarly, by repeatedly applying Lemma \ref{Precur}, we find
\begin{equation} \label{gipi2}
\theta P_{i}(\theta) = c_{i+1}c_{i+2}P_{i+2}(\theta) + 
(c_{i+1}b_{i+2}+b_{i+1}c_{i}) P_{i}(\theta) + b_{i}b_{i+1} 
P_{i-2}(\theta).  
\end{equation}
By (\ref{pigi}), we find
\begin{equation} \label{dion}
   g_{i} = 
   p_{i} + \frac{b_{i}b_{i+1}}{c_{i-1}c_{i}} 
   \frac{P_{i-2}(\theta)}{P_{i}(\theta)} g_{i-2}.
\end{equation}
Using (\ref{dion}) to eliminate $g_{i-2}$ in (\ref{gipi0}), and then applying 
(\ref{gipi1}), (\ref{gipi2}), we obtain (\ref{gipi}).
\hfill $\Box$\\

\begin{theorem} \label{P10}
With reference to Notation \ref{thetanote}  and Definition \ref{GIPOLY}, 
 for $0 \le i,j \le D-2$ we have
\begin{equation}  \label{GIGJredux}
\sum_{h=0}^{D}g_{i}(\theta_{h})g_{j}(\theta_{h})(k^{2}-\theta_{h}^{2})
(\theta -\theta_{h}^{2}) m_{h} = 
\delta_{ij}|X|k_{i}b_{i}b_{i+1}c_{i+1}c_{i+2}\frac{P_{i+2}(\theta)}
{P_{i}(\theta)}.
\end{equation}
\end{theorem}
{\it Proof.}  Without loss of generality, we may assume $i \le j$.  
First we eliminate $g_i(\theta_h)$ and $g_j(\theta_h)(\theta- \theta_h^2)$ in the left-hand side of 
(\ref{GIGJredux}) by using Definition \ref{GIPOLY} and (\ref{gipi}), respectively.  Simplifying the resulting expression using (\ref{PIPJ1}) and the fact that $i \le j$, we 
obtain the right-hand side of (\ref{GIGJredux}).  The 
result follows.  \hfill $\Box$\\

We finish this section with a comment.

\begin{lemma} \label{gD-2}
With reference to Notation \ref{thetanote}  and Definition \ref{GIPOLY}, 
assume $D$ is even and $\theta =0$.  Then $g_{D-2}(\theta_h) =0$ for $1 \leq h \leq D-1$, $\; h \not=d$.
\end{lemma}
{\it Proof.}  Recall $\theta_d =0$ by Lemma \ref{thetad}.
Setting $i=j=D-2$ and $\theta =0$ in (\ref{GIGJredux}), we find
\begin{equation} \label{i=D-2}
\sum_{{h=1}\atop{ h \not= d}}^{D-1} \theta_h^2 g_{D-2}^2(\theta_h)(k^2-\theta_h^2)m_h = -|X|k_{D-2}b_{D-2}b_{D-1}c_{D-1}c_D \frac{P_D(0)}{P_{D-2}(0)}.
\end{equation}
In (\ref{i=D-2}) the right-hand side is zero by Lemma \ref{Pnot0}.  In the left-hand side each summand is nonnegative so each summand is zero.  In each summand the factor $\theta_h^2 (k^2-\theta_h^2)m_h$ is nonzero so the remaining factor $g_{D-2}(\theta_h)$ is zero.  The result follows.
\hfill $\Box$\\


\section{The subconstituent algebra and its modules}  \label{subconsection}
In this section we recall some definitions and basic concepts 
concerning
the subconstituent algebra and its modules.
For more information we refer the reader to \cite{caugh2},
\cite{curtin1},
\cite{curtin2},
\cite{go},
\cite{hobart},
\cite{terwSub1}.

\bigskip
\noindent
Let $\Gamma=(X,R)$ denote a distance-regular graph with diameter $D 
\ge 3.$
We recall the dual Bose-Mesner algebra of $\Gamma.$
From now on we fix a vertex $x \in X.$ For $ 0 \le i 
\le D$ let $E_i^*=E_i^*(x)$ denote the diagonal
matrix in $\hbox{Mat}_X(\C)$ with $yy$ entry
\begin{equation}\label{DEFDEI}
{(E_i^*)_{yy} = \cases{1, & if $\partial(x,y)=i$\cr
0, & if $\partial(x,y) \ne i$\cr}} \qquad (y \in X).
\end{equation}
We call $E_i^*$ the {\it $i^{\hbox {th}}$ dual idempotent of} $\Gamma$
{\it with respect to x.} 
We observe
(di) $\sum_{i=0}^D E_i^* = I$;
(dii) $\overline{E_i^*} = E_i^* \; (0 \le i \le D)$; 
(diii) $E_i^{*t} = E_i^*  \;(0 \le i \le D)$;
(div) $E_i^*E_j^* = \delta_{ij}E_i^*  \;(0 \le i,j \le 
D).$
Using (di) and (div) 
we find $E_0^*,E_1^*, \ldots, E_D^*$ form a basis for a commutative 
subalgebra $M^*=M^*(x)$ of $\hbox{Mat}_X(\C).$ We call $M^*$ the {\it 
dual Bose-Mesner algebra of}
$\Gamma$ {\it with respect to x.} 
We recall the subconstituents of $\Gamma $.
Using (\ref{DEFDEI}) we find
\begin{equation}\label{DEIV}
E_i^*V = {\rm span}\,\{\hat{y} \;|\; y \in X, \quad \partial(x,y)=i\}
\qquad (0 \le i \le D).
\end{equation}
By (\ref{DEIV}) and since $\lbrace {\hat y} \;|\;y \in X\rbrace $ is
an orthonormal basis for $V$ we find 
\begin{displaymath}
V = E_0^*V+E_1^*V+ \cdots +E_D^*V \qquad \qquad {\rm (orthogonal\ 
direct\ sum}).
\end{displaymath}
Combining (\ref{DEIV}) and (\ref{DEFKI}) we find
the dimension of $E_i^*V$ is $k_i$ for $0 \le i \le D$.
We call $E_i^*V$ the {\it $i^{\hbox {th}}$ subconstituent of} $\Gamma$
{\it with respect to} $x$.

\bigskip
\noindent
We recall how $M$ and $M^*$ are related.
By \cite[Lemma 3.2]{terwSub1}, \begin{equation} \label{REL2}
E_h^*A_iE_j^*=0 \quad {\rm if\ and\ only\ if} \quad p_{ij}^h = 0
\qquad \qquad (0 \le h,i,j \le D).
\end{equation}
Combining   (\ref{REL2}) and 
 (\ref{REL1}) we find
\begin{equation}
E_i^*AE_j^*=0 \qquad {\rm if} \qquad |i-j| > 1
\qquad \qquad (0 \le i,j \le D).
\label{REL3}
\end{equation}
\noindent
Let $T=T(x)$ denote the subalgebra of $\hbox{Mat}_X(\C)$ generated by 
$M$ and $M^*$. We call $T$ the {\it subconstituent algebra of} 
$\Gamma$ {\it with respect to} $x$ \cite{terwSub1}. We observe
$T$ has finite dimension. Moreover $T$ is semi-simple;
the reason is that $T$ is closed under the conjugate-transpose map
\cite[p. 157]{CR}.

\bigskip
\noindent
We now consider the modules for $T.$ By a {\it T-module}
we mean a subspace $W \subseteq V$ such that $BW \subseteq W$
for all $B \in T.$ We refer to $V$ itself as the {\it standard
module} for $T.$ Let $W$ denote a $T$-module. Then $W$ is said
to be {\it irreducible} whenever $W$ is nonzero and $W$ contains no 
$T$-modules other than 0 and $W.$ Let $W,W^\prime$ denote
$T$-modules. By an {\it isomorphism of $T$-modules}
from $W$ to $W^\prime$ we
mean an isomorphism of vector spaces
$\sigma: W \rightarrow W^\prime$
such that
$$
(\sigma B- B \sigma)W = 0 \qquad \qquad {\rm for\ all}\; B \in T.
$$
The modules $W,W^\prime$ are said to be {\it isomorphic as 
$T$-modules}
whenever
there exists an isomorphism of $T$-modules from $W$ to $W^\prime.$

\bigskip
\noindent
Let $W$ denote a $T$-module and let $W'$ denote a $T$-module 
contained in $W$. Using (\ref{ADJ}) we find the orthogonal complement 
of $W'$ in $W$ is a $T$-module.
It follows that each $T$-module
is an orthogonal direct sum of irreducible $T$-modules.
We mention any two nonisomorphic
irreducible $T$-modules are orthogonal \cite[Chapter IV]{CR}.

\bigskip
\noindent
Let $W$ denote an irreducible $T$-module.
Using (di)--(div) above we find $W$ is the direct sum of the
nonzero spaces among $E^*_0W, E^*_1W,\ldots, E^*_DW$.
Similarly using
(eii)--(ev) we find
$W$ is the direct sum of the
nonzero spaces among $E_0W, E_1W,\ldots, E_DW$.
If the dimension of $E^*_iW$ is at most 1 for
$0 \leq i \leq D$ then
the dimension of $E_iW$ is at most 1 for
$0 \leq i \leq D$ \cite[Lemma 3.9]{terwSub1}; in this case we say $W$ 
is {\it thin}.
Let $W$ denote an irreducible $T$-module.
By the {\it endpoint}
of $W$ we mean $$ {\rm min}\,\{i\;|\;0\le i \le D, \;\; E_i^*W \ne 0 
\}.$$

\noindent
For the rest of the paper we adopt the following
notational convention.

\begin{definition}\label{A}  \rm
Let $\Gamma=(X,R)$ denote a bipartite
distance-regular graph with diameter $D \ge 4$, valency $k \ge 3$, 
intersection numbers $b_i, c_i$,
distance matrices $A_i$,
Bose-Mesner algebra $M$,
and eigenvalues
$\theta_0 > \theta_1 > \cdots > \theta_D$.
For $0 \leq i \leq D$ we let $E_i $ denote the primitive idempotent of
$\Gamma$ associated with $\theta_i$. We define $d = \lfloor D/2 \rfloor$.
We fix $x \in X$ and abbreviate
$E_i^*=E_i^*(x)$ $(0 \le i \le D)$,
$ M^*=M^*(x)$,
$T=T(x)$.  We let $V$ denote the standard 
module for $\Gamma$. 
We define
\begin{equation}
s_i=\sum_{{y \in X}\atop {\partial(x,y)=i}} {\hat y}
\qquad \qquad (0 \leq i \leq D).
\label{si}
\end{equation}
\end{definition}


\section{The $T$-module of endpoint 0}

With reference to Definition \ref{A}, there exists a unique
irreducible $T$-module with endpoint 0 \cite[Proposition 8.4]{egge1}. We 
call this module
$V_0$. The module $V_0$ is described in \cite{curtin1}, \cite{egge1}.
We summarize some details below in order to 
motivate the results that
follow.

\bigskip
\noindent The module $V_0$ is thin. In fact
each of $E_iV_0$, $E^*_iV_0$
has dimension 1 for $0 \leq i \leq D$. We give two bases for $V_0$. 
The vectors 
$E_0{\hat x}, E_1{\hat x}, \ldots , E_D{\hat x}$
form a basis for $V_0$. These vectors are mutually orthogonal and 
$\|E_i {\hat x}\|^2 = m_i|X|^{-1}$ for $0 \leq i \leq D$.
To motivate the second basis we make some comments.  
For $0 \leq i \leq D$ we have
$s_i = A_i{\hat x}$.
Moreover $s_i=E^*_i\delta $, where $\delta=\sum_{y \in X} {\hat y}$. 
The vectors 
$s_0, s_1, \ldots, s_D$
form a basis for $V_0$.
These vectors are mutually orthogonal and
$\|s_i\|^2 = k_i$ for $0 \leq i \leq D$.
With respect to the basis $s_0, s_1, \ldots, s_D$ the matrix 
representing $A$ is
\begin{eqnarray*}
\left(\begin{array}{cccccc}
0 & b_0 & & & & {\bf 0}\\
c_1 & 0 & b_1 & & & \\
& c_2 & \cdot & \cdot & & \\
& & \cdot & \cdot & \cdot& \\
& & & \cdot & \cdot & b_{D-1} \\
{\bf 0} & & & & c_{D} & 0
\end{array} \right).
\end{eqnarray*}
The two bases for $V_{0}$ given above
are related as follows.
For $0 \leq i \leq D$ we have
$$s_i = \sum_{h=0}^{D}f_i(\theta_h) E_h{\hat x},$$
where the polynomial $f_i$ is from
(\ref{FIPOLY2}).


\section{The $T$-modules of endpoint 1}

With reference to Definition \ref{A}, there exists, up to isomorphism, a unique
irreducible $T$-module with endpoint 1 \cite[Corollary 7.7]{curtin1}. 
We call this module $V_1$.
The module $V_1$ is described in \cite{curtin1}, \cite{go2}.
We summarize some details below.

\bigskip \noindent
The module $V_{1}$ is thin with dimension $D-1$.
We give two bases for $V_1$. Let $v$ denote a nonzero
vector in $E^*_1V_1$. The vectors
\begin{equation} \label{mod1basis}
E_i v \qquad (1 \leq i \leq D-1)
\end{equation}
form a basis for $V_1$ and $E_0v =0,$ $E_Dv=0$. The vectors in (\ref{mod1basis}) are mutually
orthogonal and
$$||E_iv||^2 =  \frac{m_i(k^2-\theta^2_i)}{|X|k(k-1)} ||v||^2 \qquad 
\qquad (1 \leq i \leq D-1).$$
To motivate the second basis we make some
comments.
We have
$E^*_{i+1}A_{i}v=p_i(A)v$ for $0 \leq i \leq D-1$, 
where the $p_{i}$ are from (\ref{PIPOLY}).
The vectors
\begin{equation}  \label{mod1basis2}
E^*_{i+1}A_{i}v \qquad (0 \leq i \leq D-2)
\end{equation} 
form a basis for $V_1$ and $E^*_DA_{D-1}v=0$.
The vectors in (\ref{mod1basis2}) are mutually orthogonal
and $$||E^*_{i+1}A_{i}v||^2 =
\frac{b_2 \cdots b_{i+1}}{c_1 \cdots c_{i}} ||v||^2
 \qquad \qquad (0 \leq i \leq D-2).$$
With respect to the basis (\ref{mod1basis2}) the
matrix representing $A$ is
\begin{eqnarray*}
\left(\begin{array}{cccccc}
0 & b_2 & & & & {\bf 0}\\
c_1 & 0 & b_3 & & & \\
& c_2 & \cdot & \cdot & & \\
& & \cdot & \cdot & \cdot& \\
& & & \cdot & \cdot & b_{D-1} \\
{\bf 0} & & & & c_{D-2} & 0
\end{array} \right).
\end{eqnarray*}
The two bases for $V_{1}$ given above
are related as follows.
For $0 \leq i \leq D-2$ we have
$$E^*_{i+1}A_{i}v = \sum_{h=1}^{D-1} p_i(\theta_h)E_hv.$$

\bigskip \noindent
We comment that $V_1$ appears in $V$ with multiplicity
$k-1$.
We will need the following result.

\begin{corollary} \label{V1eigval} With reference to Definition \ref{A},
let $W$ denote an irreducible $T$-module with endpoint
$1$.  Observe $E^*_2W$ is an eigenspace for $E^*_2A_2E^*_2$.
The corresponding eigenvalue is $b_3-1$.
\end{corollary}
{\it Proof.}  The desired eigenvalue is the entry in the second row 
and second column of the matrix
representing $A_2$ with respect to the basis (\ref{mod1basis2}).  
To compute this entry, first set $i=1$ 
in (\ref{AAI1}) and observe that $c_2 A_2=A^2-kI$.  Using this fact 
and the above matrix display of $A$, we verify the specified matrix 
entry is $b_{3}-1$.
\hfill $\Box$\\


\section{The local eigenvalues}

A bit later in this paper we will consider the
thin irreducible $T$-modules with endpoint 2. In order to discuss these we recall the local eigenvalues.  

\begin{definition}\label{SUBGRAPH}  \rm
With reference to Definition \ref{A},
we let $\Gamma_{2}^{2} = \Gamma_{2}^{2}(x) $ denote the
graph
$(\breve{X},\breve{R}),$ where
\begin{eqnarray}
\breve{X} &=& \{y \in X\;|\;\partial(x,y)=2\},\nonumber\\
\breve{R} &=& \{yz\;|\;y,z \in \breve{X},\, \partial(y,z)=2 \},\nonumber
\end{eqnarray}
where we recall $\partial$ denotes the path-length distance function 
for $\Gamma$.
The graph $\Gamma_{2}^{2}$ has 
exactly $k_{2}$ vertices, where $k_{2}$ is the second valency of
$\Gamma.$ Also, $\Gamma_{2}^{2}$ is regular with valency $p^{2}_{22}$.
We let $\breve {A}$ denote the adjacency matrix of $\Gamma_{2}^{2} $. The 
matrix $\breve {A}$ is symmetric with real entries; therefore $\breve 
{A}$ is diagonalizable with all eigenvalues real.
We let
$\eta_1, \eta_2 , \ldots , \eta_{k_{2}}$
denote the eigenvalues of $\breve {A}$.  We call 
$\eta_1, \eta_2 , \ldots , \eta_{k_{2}}$ the {\em local eigenvalues
of $\Gamma$ with respect to $x$}.  
\end{definition}

\noindent With reference to Definition \ref{A}, we consider the
second subconstituent $E^*_2V$. We recall 
the dimension of $E^*_2V$ is $k_{2}$.
Observe $E^*_2V$ is invariant under the action of $E^*_2A_{2}E^*_2$.
To illuminate this action we make an observation.
For an appropriate ordering of the vertices of $\Gamma$ we have $$
E^*_2A_{2}E^*_2 =
\left(\begin{array}{cc} \breve{A} & 0 \\ 0 & 0 \end{array} \right),
$$
where $\breve {A}$ is from
Definition \ref{SUBGRAPH}. Apparently the action of
$E^*_2A_{2}E^*_2$ on $E^*_2V$ is essentially the adjacency map
for $\Gamma_{2}^{2} $.
In particular the action of $E^*_2A_{2} E^*_2$ on $E^*_2V$ is 
diagonalizable
with eigenvalues $\eta_1, \eta_2, \ldots, \eta_{k_{2}}$.  We observe 
the vector $s_{2}¥$ from (\ref{si}) is contained in $E_{2}^{*}V$.  
One may easily show that $s_{2}$ is an eigenvector for $E^*_2A_{2}E^*_2$ 
with eigenvalue $p_{22}^{2}$.  
Let $v$ denote a vector in $E_{2}^{*}V$.
We observe the following are equivalent:  (i) $v$ is orthogonal to
$s_2$; (ii) $E_0v=0$; (iii) $Jv=0$; (iv) $E_{D}v=0$; (v) $J'v=0$.  
Let $V_{1}$ 
denote an irreducible $T$-module of endpoint $1$, and let $v$ 
denote a vector in $E_{2}^{*}V_{1}$.  By Corollary \ref{V1eigval}, 
$v$ is an eigenvector for $E^*_2A_{2}¥E^*_2$ with eigenvalue 
$b_{3}-1$.  
Reordering 
the local eigenvalues if necessary, we have $\eta_{1} = p_{22}^{2}$ 
and $\eta_{i}= b_{3}-1$ $(2 \leq i \leq k)$.  For the rest of this 
paper we assume the local eigenvalues of $\Gamma$ are ordered in this 
way.

\bigskip \noindent
We now need some notation.    

\begin{definition} \label{U=V0+Y}  \rm With reference to Definition \ref{A}, 
let $Y$ denote the subspace of $V$ spanned by the
   irreducible $T$-modules with endpoint 1.  We define 
$U$ to be the orthogonal 
   complement of $E^*_2V_{0} + E^*_2Y$ in $E^{*}_{2}V$.
 \end{definition}
 

\begin{definition}  \label{phidef}  \rm
With reference to Definition \ref{A}, let $\Phi$ denote the set of 
distinct scalars among $\eta_{k+1}$,$\eta_{k+2}$,$\ldots$, 
$\eta_{k_{2}}$, where the $\eta_{i}$ are from Definition 
\ref{SUBGRAPH}.  For $\eta \in \R$ we let $\mbox{mult}_{\eta}$ denote the 
number of times $\eta$ appears among 
$\eta_{k+1}, \eta_{k+2}, \ldots, \eta_{k_{2}}$.  We observe 
$\mbox{mult}_{\eta} \not= 0$ if and only if $\eta \in \Phi$.  
\end{definition}

\noindent
Using (\ref{ADJ}) we find $U$ is invariant under $E^*_2A_{2}E^*_2$. 
Apparently the restriction of $E^*_2A_{2}E^*_2$ to $U$ is diagonalizable 
with eigenvalues 
$\eta_{k+1}, \eta_{k+2}, \ldots, \eta_{k_{2}}$.
For $\eta \in \R$ let $U_\eta $ denote the set consisting of those
vectors in $U$ that are eigenvectors for $E^*_2A_{2}¥E^*_2$ with 
eigenvalue $\eta $.
We observe $U_\eta $ is a subspace of $U$ with dimension 
$\mbox{mult}_{\eta}$.
We emphasize the following are equivalent:  (i) $\mbox{mult}_{\eta} 
\not= 0$; (ii) $U_{\eta} \not= 0$; (iii) $\eta \in \Phi$.  
By (\ref{ADJ}) and since $E^*_2A_{2}¥E^*_2$ is symmetric with real entries 
we find
\begin{equation}
U = \sum_{\eta \in \Phi} U_\eta \qquad \qquad (\hbox{orthogonal 
direct sum}).
\label{Ubreakdown}
\end{equation}

\begin{definition}\label{TYPE}  \rm
With reference to Definition \ref{A},
for all $z \in \C \cup \infty$ 
we define
$${{\tilde z} = \cases{
-1-\frac{b_2 b_{3}}{z^{2}-b_{2}¥}, & if $z\not=\infty, \;z^{2} \ne {b_{2}}$ \cr
\infty, & if $z^{2} = {b_{2}}$ \cr
-1, &if $z=\infty $. \cr
}}$$
\end{definition}

\begin{note} \rm
With reference to Definition \ref{A}, neither of 
$\theta_1^{2}, \theta_d^{2}$ is equal to ${b_{2}}$ by Lemma
\ref{bipeigs},
so
\begin{equation}  \label{tildeth1d}
{\tilde \theta}_1=-1-b_2 b_{3}(\theta_1^{2}-b_{2})^{-1}, \qquad \qquad 
{\tilde \theta}_d=-1-b_2 b_{3}(\theta_d^{2}-b_{2})^{-1}. 
\end{equation}
By the data in
Lemma \ref{bipeigs} we have
${\tilde \theta}_1 < -1$. Moreover
${\tilde \theta}_d > b_3-1$ if $D$ is odd and  
${\tilde \theta}_d = b_3-1$ if $D$ is even.
In either case 
${\tilde \theta}_d \geq 0$.
\end{note}

\begin{lemma}\label{TYPEINEQ}  \cite[Theorem 11.4]{maclean5}
With reference to Definitions \ref{A}
and \ref{SUBGRAPH},
we have
${\tilde \theta}_1 \leq \eta_i \leq {\tilde \theta}_d$ for
$k+1 \le i \le k_{2}$.
\end{lemma}

\noindent We remark on the case of equality in
the above lemma.

\begin{lemma} \label{EJV}  \cite[Lemma 11.5]{maclean5}
With reference to Definition \ref{A},
let $v$ denote a nonzero vector in $U$. Then (i)--(vi) hold below:
\begin{description}
\item[{\rm (i)}] $E_0v=0$ and $E_{D}v=0$.
\item[{\rm (ii)}] For $1 \le i \le D-1$, $\; E_{i}v \not= 0$ 
provided $i$ is not among $1, d, D-d, D-1$. 
\item[{\rm (iii)}] $E_1v=0$  if and only if $v \in 
U_{{\tilde \theta}_1}$.
\item[{\rm (iv)}] $E_{D-1}v=0$ if and only if $v \in 
U_{{\tilde \theta}_1}$.
\item[{\rm (v)}] $E_dv=0$ if and only if $v \in
U_{{\tilde \theta}_d}$.
\item[{\rm (vi)}] $E_{D-d}v=0$ if and only if $v \in
U_{{\tilde \theta}_d}$.
\end{description}
\end{lemma}

\begin{corollary} \label{dimj}  \cite[Corollary 11.6]{maclean5}
With reference to Definition \ref{A},
let $v$ denote a nonzero vector in $U$.
Then (i)--(iv) hold below:
\begin{description}
\item[{\rm (i)}] If $v \in
U_{{\tilde \theta}_1}$ then
$Mv$ has dimension $D-3$.
\item[{\rm (ii)}] If $v \in
U_{{\tilde \theta}_d}$ and $D$ is odd, then
$Mv$ has dimension $D-3$.
\item[{\rm (iii)}] If $v \in
U_{{\tilde \theta}_d}$ and $D$ is even, then
$Mv$ has dimension $D-2$.
\item[{\rm (iv)}]
If $v \notin
U_{{\tilde \theta}_1}$
and
$v \notin
U_{{\tilde \theta}_d}$ then
$Mv$ has dimension $D-1$.
\end{description}
\end{corollary}

\begin{definition}
\label{TM}  \rm
With reference to Definition \ref{A}, let
$W$ denote a thin irreducible $T$-module with endpoint 2.
Observe $E^*_2W$ is a $1$-dimensional eigenspace for
$E^*_2A_{2}¥E^*_2$; let $\eta $ denote the corresponding eigenvalue.
We observe $E^*_2W$ is contained in $E^*_2V$ and is orthogonal
to any irreducible $T$-module with endpoint 0 or 1, so
$E^*_2W\subseteq U_\eta$.
Apparently $U_\eta \not=0$ so $\eta $
is among $\eta_{k+1}, \eta_{k+2},\ldots, \eta_{k_{2}}$.
We have ${\tilde \theta}_1 \leq
\eta \leq {\tilde \theta}_d$ by Lemma \ref{TYPEINEQ}. We refer to 
$\eta $ as the {\em local eigenvalue} of $W$.
\end{definition}

\noindent With reference to Definition \ref{A}, let $W$ denote a thin 
irreducible $T$-module with endpoint 2 and local eigenvalue
$\eta $. In order to describe $W$ we distinguish four cases: (i) $\eta 
={\tilde \theta}_1 $; (ii) $D$ is odd 
and $\eta = {\tilde \theta}_{d}$; (iii) $D$ is even 
and $\eta = {\tilde \theta}_{d}$; (iv) $ {\tilde \theta}_1 < \eta < {\tilde 
\theta}_d$.  
For cases (i), (ii) the module $W$ was described by the present 
authors in \cite{maclean5}; we summarize these results in the 
following section.  For cases (iii), (iv) we describe $W$ in Sections  \ref{main2} and \ref{main1}.


\section{Some thin irreducible $T$-modules with endpoint 2}   

\noindent 
In this section we summarize some results from
\cite{maclean5}
concerning the thin irreducible $T$-modules with endpoint 2 and
local eigenvalue $\eta$, where 
$\eta = {\tilde \theta}_1$,  or 
$\eta ={\tilde \theta}_d$ with $D$ odd.

\medskip
\noindent
With reference to Definition \ref{A}, choose 
$n \in \{1,d\}$ if $D$ is odd, and let $n=1$ if $D$ is even.  
Define $\eta = {\tilde \theta}_{n}$.  
Let $W$ denote a thin irreducible
$T$-module with endpoint 2 and local eigenvalue $\eta$.  
The dimension of
$W$ is $D-3$. For $0 \leq i \leq D$, 
 $E^*_iW$ is zero if 
$i\in \lbrace 0,1,D-1,D\rbrace $, and has dimension 1 if 
$i \not\in \lbrace 0,1,D-1,D\rbrace $. 
Moreover 
$E_iW$ is zero if
$i\in \lbrace 0,n,D-n,D\rbrace $, and has dimension 1 if 
$i \not\in \lbrace 0,n,D-n,D\rbrace $.
Let $v$ denote
a nonzero vector in $E^*_2W$. Then $W=Mv$.  
The vectors  
\begin{equation}\label{short5a}
E_iv  \qquad  \qquad (1 \leq i \leq D-1, \;\; i\not=n, \; i \not= D-n)
\end{equation}
form a basis for $W$, and each of $E_0v, E_nv, E_{D-n}v, E_{D}v$ is 
zero. 
The vectors in  
(\ref{short5a})  
are mutually orthogonal and
$$
\|E_iv\|^2 = \frac{m_i(\theta_i^{2}-k^{2})(\theta_i^{2}-\theta_{n}^{2})}
{|X|k b_{1}(\theta_{n}^{2} - b_{2})}
\|v\|^2
\qquad \qquad (1 \le i \le D-1, \;\;i\not=n, \;\; i\not= D-n).
$$
We mention a second basis for $W$. To motivate things we remark
$$
E^*_{i+2}A_iv = \sum_{{h=0}\atop{ i-h { \mbox{ \tiny even}}}}^{i}
\frac{p_h(\theta_{n})}{p_i(\theta_n)}\frac{k_i b_i b_{i+1}}{ k_hb_h b_{h+1}} 
p_h(A)v \qquad (0 \le i \le D-2).
$$
The vectors  
\begin{equation}\label{short2a}
E^*_{i+2}A_iv \qquad \qquad (0 \leq i \leq D-4)
\end{equation}
form a basis for $W$, and $E^*_{D-1}A_{D-3}v=0, \; E^*_DA_{D-2}v=0$.  
The vectors in (\ref{short2a}) are mutually  orthogonal and
$$
\|E^*_{i+2}A_iv\|^2 = 
\frac{k_{i}b_{i}b_{i+1}c_{i+1}c_{i+2}}{kb_{1}(\theta_{n}^{2} - b_{2})}
\frac{p_{i+2}(\theta_{n})}{p_{i}(\theta_{n})} \| v \|^{2}  \qquad (0 
\le i \le D-4).
$$
With respect to the basis given in 
(\ref{short2a}) the matrix representing $A$ is
$$
\left(\begin{array}{cccccc}
0 & w_1 & & & & {\bf 0}\\
c_1 & 0 & w_2 & & & \\
& c_2 & \cdot & \cdot & & \\
& & \cdot & \cdot & \cdot& \\
& & & \cdot & \cdot & w_{D-4} \\
{\bf 0} & & & & c_{D-4} & 0
\end{array} \right),
$$
where 
\begin{displaymath}  
  w_{i} = \frac{b_{i+1}c_{i+2}}{c_{i}} \frac{p_{i-1}(\theta_{n}) 
  p_{i+2}(\theta_{n})}{p_{i}(\theta_{n}) p_{i+1}(\theta_{n})} \qquad 
  \qquad (1 \le i \le D-4).
\end{displaymath}
The bases for $W$ given in 
(\ref{short5a}), 
(\ref{short2a})
are related as follows.
For $0 \leq i \leq D-4$ we have 
$$
E^*_{i+2}A_iv = \sum_{{1 \le j \le D-1}\atop{j \not=n, \, j \not=D-n}}
\gamma_i(\theta_j) 
E_jv,
$$
where
$$
\gamma_i = \sum_{{h=0}\atop{ i-h { \mbox{ \tiny even}}}}^{i}
\frac{p_h(\theta_{n})}{p_i(\theta_n)}\frac{k_i b_i b_{i+1}}{ k_hb_h b_{h+1}} p_h.
$$

\noindent 
We finish this section with a comment.

\begin{lemma}\cite[Theorem 12.9]{maclean5}
\label{lem:mvisshort}
With reference to Definition \ref{A}, let 
$v$ denote a nonzero vector in $U$.  
Let
$n \in \{1,d\}$ if $D$ is odd, and let $n=1$ if $D$ is even.
Assume $v$ is an eigenvector for $E^*_2A_{2}E^*_2$
with eigenvalue
${\tilde \theta}_n$.
 Then $Mv$ is a thin irreducible $T$-module
with endpoint 2 and local eigenvalue 
${\tilde \theta}_n$.
\end{lemma}

\section{The space $Mv$ when $D$ is even and $v \in U_{{\tilde \theta}_d}$}

With reference to Definition \ref{A}, assume $D$ is even.
One of our ultimate goals in this paper is to describe 
the thin irreducible $T$-modules with endpoint 2 and 
local eigenvalue ${\tilde \theta}_d$. 
Before we get to this, we find it illuminating to 
 consider a more general type of space. 
Let $v$ denote a nonzero vector
in $U$ and assume $v$ is an eigenvector for $E^*_2A_{2}E^*_2$
with corresponding eigenvalue ${\tilde \theta}_d$. 
In this section we investigate the
space $Mv$.
We present two orthogonal bases for $Mv$ which we find attractive.
Recall that since $D$ is even, we have $\theta_d =0$ and thus ${\tilde \theta}_d = b_3-1$.

\begin{theorem} 
\label{M5d}
With reference to Definition \ref{A}, assume $D$ is even, and let
$v$ denote a nonzero vector in $U$. Assume $v$ is an eigenvector
for $E^*_2A_{2}E^*_2$ with corresponding
eigenvalue ${\tilde \theta}_d$.  Then the vectors  
$E_iv$ $(1 \le i \le D-1, \; i \not= d)$
form a basis for $Mv$. Moreover $E_0v=0$, $E_d v =0$, $E_{D}v=0$.
\end{theorem}
{\it Proof.} Recall $E_0, E_1, \ldots, E_D$ form a basis
for $M$. 
 Observe  $E_0v=0$, $E_d v =0$, $E_{D}v=0$ 
by 
 Lemma \ref{EJV}
so the vectors
$E_iv$  $(1 \le i \le D-1, \; i \not= d)$ span $Mv$.
 These vectors are nonzero 
by 
Lemma \ref{EJV} and mutually orthogonal by
(\ref{EIVDECOM}), so they are linearly independent.
The result follows.
\hfill $\Box$\\

\begin{theorem}  \cite[Theorem 11.2]{maclean5}
\label{M6d}
With reference to Definition \ref{A}, assume $D$ is even, and let
$v$ denote a nonzero vector in $U$. Assume $v$ is an eigenvector
for $E^*_2A_{2}E^*_2$  with corresponding
eigenvalue ${\tilde \theta}_d$.  Then the vectors 
$E_iv$ $(1 \le i \le D-1, \; i \not= d)$ are mutually orthogonal.
Moreover the square-norms of these vectors are given
as follows:
 \begin{displaymath} 
\|E_iv\|^2
= \frac{m_i(k -\theta_i)(k + \theta_{i})\theta_i^{2}}
{|X|k b_{1}b_{2}}
\|v\|^2
\qquad \qquad (1 \le i \le D-1, \; i \not= d).
\end{displaymath} 
(The scalar $m_i$ denotes the multiplicity of $\theta_i$.)
\end{theorem}

\noindent Referring to Theorem \ref{M5d}, 
we now consider a second basis for $Mv$.

\begin{definition} 
\label{M1d}  \rm
With reference to Definition \ref{A}, assume $D$ is even, and let
$v$ denote a nonzero vector in $U$. Assume $v$ is an eigenvector
for $E^*_2A_{2}E^*_2$ with corresponding
eigenvalue ${\tilde \theta}_d$. 
We define the vectors $v_0, v_1, \ldots, v_{D-2}$ by
\begin{equation}\label{M1ad}
v_i = \sum^{i}_{{h=0}\atop{ i-h { \mbox{ \tiny even}}}}
\frac{P_h(0)}{P_i(0)}\,
\frac{k_ib_ib_{i+1}}{k_hb_hb_{h+1}} \,p_h(A)v
\qquad \qquad (0 \leq i \leq D-2).
\end{equation}
(The polynomials $p_i$ are from  
(\ref{PIPOLY}), and the $P_{i}$ are from (\ref{Pidef2}).)
The denominators in (\ref{M1ad}) are nonzero by 
Corollary \ref{Pnot0v2}.
\end{definition}

\begin{theorem} 
\label{M2d}
With reference to Definition \ref{A}, assume $D$ is even, and let
$v$ denote a nonzero vector in $U$. Assume $v$ is an eigenvector
for $E^*_2A_{2}E^*_2$ with corresponding
eigenvalue ${\tilde \theta}_d$. Then with reference to (\ref{M1ad}), the vectors 
$v_0, v_1, \ldots , v_{D-3}$ form a basis for $Mv$ and $v_{D-2}=0$.
\end{theorem}
{\it Proof.} 
By  Theorem
\ref{M5d}
  we find $Mv$ has dimension $D-2$. By this
and since $A$ generates $M$, we find $Mv$ has a  basis
$v, Av, \ldots, A^{D-3}v$.
For $0 \leq i \leq D-3$ the vector  $v_i$ is contained in the span of
$v, Av, \ldots,  A^iv$ but not in the span of
$v, Av, \ldots,  A^{i-1}v$.
It follows that $v_0, v_1, \ldots, v_{D-3}$ form a basis for $Mv$.
To see that $v_{D-2}=0$, first let $g_{D-2}$ denote the polynomial from 
Definition \ref{GIPOLY}, where $\theta=0$.  Comparing (\ref{GIPOLY1}),
(\ref{M1ad}) we find $v_{D-2}=g_{D-2}(A)v$. Using this and 
(eii) we routinely obtain
$v_{D-2} =
\sum_{j=0}^D g_{D-2}(\theta_j) E_jv$.  Applying Lemma \ref{gD-2} and 
Theorem \ref{M5d}, we find $v_{D-2}=0$.
\hfill $\Box$\\

\noindent With reference to
Definition
\ref{M1d},
we will show the vectors 
$v_0, v_1, \ldots , v_{D-3}$
are mutually orthogonal and we will compute their 
square-norms. To do this we need the following result. 

\begin{theorem} 
\label{M7d}
With reference to Definition \ref{A}, assume $D$ is even, and let
$v$ denote a nonzero vector in $U$. Assume $v$ is an eigenvector
for $E^*_2A_{2}E^*_2$ with corresponding
eigenvalue ${\tilde \theta}_d$. 
Let the vectors $v_0, v_1, \ldots, v_{D-3}$
be as in Definition
\ref{M1d}.
Then for $0 \leq i \leq D-3$ we have 
\begin{equation}\label{M7ad}
v_i =   \sum_{{j=1}\atop{j \not= d}}^{D-1} g_i(\theta_j) E_jv,
\end{equation}
where
\begin{equation}\label{M7bd}
g_i = \sum_{{h=0}\atop{ i-h { \mbox{ \tiny even}}}}^{i} 
\frac{P_h(0)}{P_i(0)}\,
\frac{k_ib_ib_{i+1}}{k_hb_hb_{h+1}} \,p_h.
\end{equation}
\end{theorem}
{\it Proof.} Let the integer $i$ be given. Comparing 
(\ref{M1ad}),
(\ref{M7bd}) we find $v_i=g_i(A)v$. Using this and 
(eii) we routinely obtain
$v_i =
\sum_{j=0}^D g_i(\theta_j) E_jv$.
Line  (\ref{M7ad}) follows since $E_0v=0$, $E_dv=0$, $E_Dv=0$ by Theorem \ref{M5d}. 
\hfill $\Box$\\

\begin{theorem} 
\label{M3d}
With reference to Definition \ref{A}, assume $D$ is even, and let
$v$ denote a nonzero vector in $U$. Assume $v$ is an eigenvector
for $E^*_2A_{2}E^*_2$ with corresponding
eigenvalue ${\tilde \theta}_d$. 
Then the vectors
 $v_0, v_1, \ldots, v_{D-3}$
from Definition
\ref{M1d}
 are mutually orthogonal.
Moreover the square-norms of these vectors are given
as follows:
\begin{equation}\label{M3ad}
\|v_i\|^2 = - \frac{k_{i}b_{i}b_{i+1}c_{i+1}c_{i+2}}{k b_{1}b_{2}} \,
\frac{P_{i+2}(0)}
{P_i(0)} 
\, \|v \|^2
\qquad \qquad (0 \leq i \leq D-3).
\end{equation}
\end{theorem}
{\it Proof.} 
Let the 
polynomials $g_0, g_1, \ldots, g_{D-3}$ be as in
(\ref{M7bd}). Using in order
Theorem 
\ref{M7d},
Theorem \ref{M6d}, 
and 
Theorem \ref{P10},
we find that for $0 \leq i,j\leq D-3$, 
\begin{eqnarray}
\langle v_i, v_j \rangle &=&  
\sum_{{h=1}\atop{h \not= d}}^{D-1} g_i(\theta_h)g_j(\theta_h)\|E_hv\|^2
\nonumber \\
&=& \sum_{{h=1}\atop{h \not= d}}^{D-1} g_i(\theta_h)g_j(\theta_h)\,
\frac{m_h(k -\theta_h 
)(k + \theta_{h})\theta_h^{2}}{|X|kb_{1}b_{2}} \|v\|^2
\nonumber \\
&=& - \delta_{ij} \, \frac{k_{i}b_{i}b_{i+1}c_{i+1}c_{i+2}}{kb_{1}b_{2}} \,
\frac{P_{i+2}(0)}
{P_i(0)}
\, \|v \|^2.
\nonumber
\end{eqnarray}
Apparently $v_0, v_1, \ldots, v_{D-3}$ are mutually orthogonal
and satisfy 
(\ref{M3ad}). 
\hfill $\Box$\\

\begin{theorem} 
\label{M4dd}
With reference to Definition \ref{A}, assume $D$ is even, and let
$v$ denote a nonzero vector in $U$. Assume $v$ is an eigenvector
for $E^*_2A_{2}E^*_2$ with corresponding
eigenvalue ${\tilde \theta}_d$. 
With respect to the basis 
$v_0, v_1, \ldots, v_{D-3}$
for $Mv$ given in Definition
\ref{M1d}
the matrix representing $A$ is
$$
\left(\begin{array}{cccccc}
0 & \omega_1 & & & & {\bf 0}\\
c_1 & 0 & \omega_2 & & & \\
 & c_2 & \cdot & \cdot & & \\
 & & \cdot & \cdot &  \cdot& \\
 & & & \cdot & \cdot & \omega_{D-3} \\
 {\bf 0} & & & & c_{D-3} & 0
\end{array} \right),
$$
where 
\begin{equation}
\omega_i = 
\label{M4aad}  
\,\frac{ b_{i+1}c_{i+2}}{c_i}
\frac{P_{i-1}(0)P_{i+2}(0)}
{P_{i}(0) P_{i+1}(0)}
\qquad \qquad (1 \leq i \leq D-3).
\end{equation}
\end{theorem}
{\it Proof.} For  $0 \leq i \leq D-2$ we define $g_i$ as in
Definition \ref{GIPOLY}, where $\theta =0$. 
 Setting 
$\lambda =A$ and $\theta =0$ in Theorem
\ref{P9} we find
\begin{equation}\label{M4dv2}
Ag_i(A) = c_{i+1}g_{i+1}(A) +  \omega_ig_{i-1}(A) 
\qquad \qquad (0 \le i \le D-3),
\end{equation}
where $g_{-1}=0$, $\omega_{0}=0$, and the $\omega_i$ are
from 
(\ref{M4aad}).
 Observe $g_i(A)v=v_i$ for $0 \leq i \leq D-2$. 
 Applying (\ref{M4dv2}) to $v$, and simplifying the result using
these comments, we find
$$
Av_i = c_{i+1}v_{i+1} + \omega_iv_{i-1} 
\qquad \qquad (0 \le i \le D-3),
$$
where $v_{-1}=0$. 
The result follows from this and since $v_{D-2}=0$ by Theorem \ref{M2d}.
\hfill $\Box$\\


\section{The thin irreducible $T$-modules with endpoint 2
and local eigenvalue $ {\tilde \theta}_d$, when $D$ is even}  \label{main2}

\noindent 
With reference to Definition \ref{A}, assume $D$ is even.  We now
describe the thin irreducible $T$-modules
with endpoint 2
and local eigenvalue $ {\tilde \theta}_d$.   This section contains some of our main results.  Because of this we have tried to make it as self-contained as possible.

\begin{theorem} 
\label{T5d}
With reference to Definition \ref{A}, assume $D$ is even, and let
$W$ denote a thin irreducible $T$-module with  endpoint 2
and local eigenvalue ${\tilde \theta}_d$. 
Let $v$ denote
a nonzero vector in $E^*_2W$. 
Then $W=Mv$. The vectors  
\begin{equation}\label{T5ad}
E_iv \qquad (1 \le i \le D-1, \; i \not= d)
\end{equation}
form a basis for $W$ and $E_0v=0$, $E_d v =0$, $E_{D}v=0$.
\end{theorem}
{\it Proof.} 
We first show $W=Mv$.
From the construction $Mv$ is nonzero and contained in $W$.
Consequently in order to show $Mv=W$, it suffices
to show $Mv$ is a $T$-module.
By construction $Mv$ is closed under multiplication by $M$.
We now show that $Mv$ is closed under multiplication by $M^*$.
By Definition \ref{TM} the vector $v$ is contained in $U$. 
Moreover $v$ is an eigenvector
for $E^*_2A_{2}E^*_2$ with eigenvalue ${\tilde \theta}_d$.
Observe that $Mv$ has basis 
$v, Av, \ldots, A^{D-3}v$ by Definition \ref{M1d} and
Theorem \ref{M2d}.
Using this and 
(\ref{REL3})
we find 
$Mv \subseteq \sum_{h=2}^{D-1} E^*_hW$.
Observe the dimension of $Mv$ is $D-2$ and the dimension
of  $\sum_{h=2}^{D-1} E^*_hW$ is at most $D-2$.
Therefore 
$Mv = \sum_{h=2}^{D-1} E^*_hW$.
From this we find
$Mv$ is closed under multiplication by $M^*$ as desired.
We have shown that $Mv$ is a nonzero $T$-submodule of $W$ 
so $Mv=W$ by the irreducibility of $W$.
The remaining assertions of the present theorem follow in view of Theorem \ref{M5d}.
 \hfill $\Box$\\

\begin{theorem} 
\label{T6d}
With reference to Definition \ref{A}, assume $D$ is even, and let
$W$ denote a thin irreducible $T$-module with  endpoint 2
and local eigenvalue $ {\tilde \theta}_d$. 
Then the basis vectors for $W$ from  
(\ref{T5ad})  
are mutually orthogonal.
Moreover the square-norms of these vectors are given
as follows:
$$\|E_iv\|^2
= \frac{m_i(k - \theta_i)(k + \theta_{i})\theta_i^{2}}
{|X|k b_{1}b_{2}}
\|v\|^2
\qquad \qquad (1 \le i \le D-1, \; i \not= d).$$
(The scalar $m_i$ denotes the multiplicity of $\theta_i$.)
\end{theorem}
{\it Proof.}  
By Definition \ref{TM} the vector $v$ is contained in $U$. 
Moreover $v$ is an eigenvector
for $E^*_2A_{2}E^*_2$ with eigenvalue ${\tilde \theta}_d$.
Applying  Theorem \ref{M6d} we obtain the result.
\hfill $\Box$\\

\begin{theorem} 
\label{T1d}
With reference to Definition \ref{A}, assume $D$ is even, and let
$W$ denote a thin irreducible $T$-module with  endpoint 2
 and local eigenvalue $ {\tilde \theta}_d$. 
 Let $v$ denote
a nonzero vector in $E^*_2W$.  Then
\begin{equation} \label{celine}
E^*_{i+2}A_iv = \sum_{{h=0}\atop{ i-h { \mbox{ \tiny even}}}}^{i} 
\frac{P_h(0)}{P_i(0)}\,
\frac{k_ib_ib_{i+1}}{k_hb_hb_{h+1}} \,p_h(A)v
\qquad \qquad (0 \leq i \leq D-2).
\end{equation}
Moreover, each side of (\ref{celine}) is zero for $i=D-2$.
(The polynomials $p_i$ are from 
(\ref{PIPOLY}), and the $P_{i}$ are from (\ref{Pidef2}).)
\end{theorem}
{\it Proof.}
By Definition \ref{TM} the vector $v$ is contained in $U$. 
Moreover $v$ is an eigenvector
for $E^*_2A_{2}E^*_2$ with eigenvalue ${\tilde \theta}_d$.
Let the vectors $v_0, v_1, \ldots, v_{D-2}$
be as in Definition 
\ref{M1d}. 
We show $E^*_{i+2}A_iv=v_i$ for $0 \leq i \leq D-2$.
Using 
(\ref{REL3})
we find 
$A^iv$ is contained in $E^*_2W+\cdots + E^*_{i+2}W$
for $0 \leq i \leq D-2$. Also for 
 $0 \leq i \leq D-2$, 
 $v_i$ is a linear combination
of $v, Av, \ldots, A^iv $, so $v_i $ is contained in 
$E^*_2W+\cdots + E^*_{i+2}W$.
By this and since  $v_0, v_1, \ldots, v_{D-3}$
are  linearly independent, we find 
\begin{equation}\label{mpd}
v_0, v_1, \ldots, v_i \quad \hbox{is a basis for }   \quad 
E^*_2W+E^*_3W+\cdots + E^*_{i+2}W \qquad  (0 \leq i \leq D-3).
\end{equation}
For the rest of this proof, 
fix an integer $i$ $(0 \leq i \leq D-2)$. 
We show $v_i $ is contained in $E^*_{i+2}W$.
To see this,
recall $E^*_2W,\ldots, E^*_{D}W$
are mutually orthogonal. Therefore $E^*_{i+2}W $ is equal to 
the orthogonal complement of 
$E^*_2W+\cdots + E^*_{i+1}W $ in 
$E^*_2W+\cdots + E^*_{i+2}W$. Recall $v_i$ is orthogonal to 
each of $v_0, v_1, \ldots, v_{i-1}$. By (\ref{mpd}) the
vectors 
 $v_0, v_1, \ldots, v_{i-1}$ form a basis for  
$E^*_2W+\cdots + E^*_{i+1}W $ so $v_i$ is orthogonal to  
$E^*_2W+\cdots + E^*_{i+1}W $. Apparently $v_i$ is 
contained in $E^*_{i+2}W$ as desired.
We show $E^*_{i+2}A_iv=v_i$.
We mentioned 
the vector $v_i$ is a linear combination 
of $v, Av, \ldots, A^iv$. In this combination 
the coefficient of $A^iv$ is $(c_1 c_2\cdots c_i)^{-1}$ 
in view of Lemma \ref{P8}(ii). Similarly $A_iv$ 
is a linear combination of
$v, Av, \ldots, A^iv$, and in this combination   
the coefficient of $A^iv$ is $(c_1 c_2\cdots c_i)^{-1}$. 
Apparently $A_iv-v_i$ is a linear combination of 
$v, Av, \ldots, A^{i-1}v$. From this and our above comments
$A_iv-v_i$ is contained in  
$E^*_2W+\cdots + E^*_{i+1}W $ so $E^*_{i+2}(A_iv-v_i)$ is
zero. We already  showed $v_i \in E^*_{i+2}W$ so
 $E^*_{i+2}v_i=v_i$. Now
$E^*_{i+2}A_iv=v_i$ as desired.   Recall $v_{D-2}=0$ by Theorem \ref{M2d}, so both sides
of (\ref{celine}) are zero for $i=D-2$. 
\hfill $\Box$\\

\begin{theorem} 
\label{T2d}
With reference to Definition \ref{A}, assume $D$ is even, and let
$W$ denote a thin irreducible $T$-module with  endpoint 2
 and local eigenvalue ${\tilde \theta}_d$. 
 Let $v$ denote
a nonzero vector in $E^*_2W$. 
Then the vectors  
\begin{equation}\label{T2ad}
E^*_{i+2}A_iv \qquad \qquad (0 \leq i \leq D-3)
\end{equation}
form a basis for $W$.
\end{theorem}
{\it Proof.}  By Definition 
\ref{TM} the vector $v$ is contained in $U$. 
Moreover $v$ is an eigenvector for $E^*_2A_{2}E^*_2$
with eigenvalue ${\tilde \theta}_d$. Let the vectors $v_0, v_1, \ldots, v_{D-3}$
be as in Definition
\ref{M1d}.  By 
 Theorem
\ref{M2d} the vectors 
$v_0, v_1, \ldots, v_{D-3}$
form a basis for $Mv$. Recall $Mv=W$ by Theorem
\ref{T5d} so 
$v_0, v_1, \ldots, v_{D-3}$
form a basis for $W$. By
Theorem \ref{T1d}
$v_i = E^*_{i+2}A_iv$ for $0 \leq i \leq D-3$ and
the result follows.
\hfill $\Box$\\

\begin{theorem} 
\label{T3d}
With reference to Definition \ref{A}, assume $D$ is even, and let
$W$ denote a thin irreducible $T$-module with  endpoint 2
and local eigenvalue ${\tilde \theta}_d$. 
Then the vectors in (\ref{T2ad}) are mutually  orthogonal.
Moreover the square-norms of these vectors are given
as follows:
$$
\|E^*_{i+2}A_iv\|^2 =   -\frac{k_{i}b_{i}b_{i+1}c_{i+1}c_{i+2}}
{k b_{1}b_{2}} \,
\frac{P_{i+2}(0)}
{P_i(0)} 
\, \|v \|^2
\qquad \qquad (0 \leq i \leq D-3).
$$
\end{theorem}
{\it Proof.}
  By Definition 
\ref{TM} the vector $v$ is contained in $U$. 
Moreover $v$ is an eigenvector for $E^*_2A_{2}E^*_2$
with eigenvalue ${\tilde \theta}_d$. 
The result follows in view of 
Theorem  
\ref{M3d} and Theorem
\ref{T1d}.
\hfill $\Box$\\

\begin{theorem} 
\label{T4d}
With reference to Definition \ref{A}, assume $D$ is even, and let
$W$ denote a thin irreducible $T$-module with  endpoint 2
 and local eigenvalue ${\tilde \theta}_d$. 
With respect to the basis for $W$ given in 
(\ref{T2ad}) the matrix representing $A$ is
$$
\left(\begin{array}{cccccc}
0 & \omega_1 & & & & {\bf 0}\\
c_1 & 0 & \omega_2 & & & \\
 & c_2 & \cdot & \cdot & & \\
 & & \cdot & \cdot &  \cdot& \\
 & & & \cdot & \cdot & \omega_{D-3} \\
 {\bf 0} & & & & c_{D-3} & 0
\end{array} \right),
$$
where 
\begin{equation} \nonumber
\omega_i =  
\,\frac{ b_{i+1}c_{i+2}}{c_i}
\frac{P_{i-1}(0)P_{i+2}(0)}
{P_{i}(0) P_{i+1}(0)}
\qquad \qquad (1 \leq i \leq D-3).
\end{equation}
\end{theorem}
{\it Proof.} 
  By Definition 
\ref{TM} the vector $v$ is contained in $U$. 
Moreover $v$ is an eigenvector for $E^*_2A_{2}E^*_2$
with eigenvalue ${\tilde \theta}_d$. 
The result follows in view of 
Theorem  
\ref{M4dd} 
and Theorem
\ref{T1d}.
\hfill $\Box$\\

\begin{theorem} 
\label{T7d}
With reference to Definition \ref{A}, assume $D$ is even, and let
$W$ denote a thin irreducible $T$-module with  endpoint 2
 and local eigenvalue ${\tilde \theta}_d$. 
Let $v$ denote a nonzero vector in $E^*_2W$.
Then for $0 \leq i \leq D-3$ we have 
$$
E^*_{i+2}A_iv =   \sum_{{j=1}\atop{j \not= d}}^{D-1} g_i(\theta_j) E_jv,
$$
where
$$
g_i = \sum_{{h=0}\atop{ i-h { \mbox{ \tiny even}}}}^{i} 
\frac{P_h(0)}{P_i(0)}\,
\frac{k_ib_ib_{i+1}}{k_hb_hb_{h+1}} \,p_h.
$$
\end{theorem}
{\it Proof.} 
  By Definition 
\ref{TM} the vector $v$ is contained in $U$. 
Moreover $v$ is an eigenvector for $E^*_2A_{2}E^*_2$
with eigenvalue ${\tilde \theta}_d$. 
The result follows in view of 
 Theorem \ref{M7d} and  
Theorem \ref{T1d}.
\hfill $\Box$\\

\noindent In summary we have the following theorem.

\begin{theorem}
\label{SUMd}
With reference to Definition \ref{A}, assume $D$ is even, and let
$W$ denote a thin irreducible $T$-module with  endpoint 2
 and local eigenvalue ${\tilde \theta}_d$. 
Then $W$ has dimension $D-2$.
For $0 \leq i \leq D$, $E^*_iW$ is zero if $i \in \{0,1,D \}$ and has dimension
1 if $2 \leq i \leq D-1$. Moreover  
 $E_iW$ is zero if $i \in \{0, d, D \}$ and has dimension
1 if $1 \leq i \leq D-1, \; i \not= d$. 
\end{theorem}
{\it Proof.} The dimension of $W$ is 
$D-2$ by Theorem
\ref{T5d}.
Fix an integer $i$ $(0 \leq i \leq D)$.
From Theorem
\ref{T2d} we find 
 $E^*_iW$ is zero if $i \in \{0, 1, D \}$ and has dimension
1 if $2 \leq i \leq D-1$.   
From Theorem
\ref{T5d}
we find 
 $E_iW$ is zero if $i \in \{0, d, D \}$ and has dimension
1 if $1 \leq i \leq D-1, \; i \not= d$.   
\hfill $\Box$\\

\section{The space $Mv$ for $v \in U_{\eta} \quad ( {\tilde \theta}_1 <  
\eta < {\tilde \theta}_d)$}

\noindent 
With reference to Definition \ref{A}, let $v$ denote a nonzero vector in $U$. Assume $v$ is an eigenvector
for $E^*_2A_{2}E^*_2$, and let $\eta$ denote the corresponding
eigenvalue.  Assume
$ {\tilde \theta}_1 <  
\eta < {\tilde \theta}_d$.  Given these assumptions we will examine the space $Mv$.

\begin{theorem} 
\label{M5}
With reference to Definition \ref{A}, let
$v$ denote a nonzero vector in $U$. Assume $v$ is an eigenvector
for $E^*_2A_{2}E^*_2$ and let $\eta$ denote the corresponding
eigenvalue. Assume 
$ {\tilde \theta}_1 <  
\eta < {\tilde \theta}_d$. 
 Then the vectors  
$E_1v, E_2v, \ldots, E_{D-1}v$ 
form a basis for $Mv$. Moreover $E_0v=0$, $E_{D}v=0$.
\end{theorem}
{\it Proof.} Recall $E_0, E_1, \ldots, E_D$ form a basis
for $M$. 
 Observe  $E_0v=0$, $E_{D}v=0$ 
by 
 Lemma \ref{EJV}
so 
$E_1v, E_2v, \ldots, E_{D-1}v$  span $Mv$.
 These vectors are nonzero 
by 
Lemma \ref{EJV} and mutually orthogonal by
(\ref{EIVDECOM}), so they are linearly independent.
The result follows.
\hfill $\Box$\\

\begin{theorem}  \cite[Theorem 11.2]{maclean5}
\label{M6}
With reference to Definition \ref{A}, let
$v$ denote a nonzero vector in $U$. Assume $v$ is an eigenvector
for $E^*_2A_{2}E^*_2$ and let $\eta$ denote the corresponding
eigenvalue. Assume
$ {\tilde \theta}_1 <  
\eta < {\tilde \theta}_d$. 
Then the vectors 
$E_1v, E_2v, \ldots, E_{D-1}v$ 
are mutually orthogonal.
Moreover the square-norms of these vectors are given
as follows:
\begin{description}
\item[{\rm (i)}]
Assume $\eta \not=-1$. Then \begin{equation}\label{EIV1}
\|E_iv\|^2
= \frac{m_i(\theta_i-k)(\theta_{i}+k)(\theta_i^{2}-\psi)}
{|X|k b_{1}(\psi - b_{2})}
\|v\|^2
\qquad \qquad (1 \le i \le D-1),
\end{equation} where \begin{equation}  \label{psidef} \psi = b_{2}\left(1- 
\frac{b_{3}}{1+\eta}\right). \end{equation}
We remark the denominator in (\ref{EIV1}) is nonzero by
(\ref{psidef}).  
\item[{\rm (ii)}]
Assume $\eta =-1$. Then \begin{displaymath} 
\|E_iv\|^2
= \frac{m_i(k - \theta_i)(k + \theta_i)}{|X|kb_1} \|v\|^2
\qquad \qquad (1 \le i \le D-1).
\end{displaymath}
\end{description}
(The scalar $m_i$ denotes the multiplicity of $\theta_i$.)
\end{theorem}

\noindent As we proceed in this section, we will encounter scalars
of the form $P_i(\psi)$ in the denominator of
some rational expressions. To make it clear these scalars are nonzero
we present the following result. 

\begin{lemma}
\label{M0}
With reference to Definition \ref{A}, let
$\eta $ denote a real number such that $\eta \not= -1$ and
${\tilde \theta}_1 < 
\eta < {\tilde \theta}_d$, and let $\psi$ be as in (\ref{psidef}).
Then (i)--(iii) hold below:
\begin{description}
\item[{\rm (i)}]  
Assume ${\tilde \theta}_1 < 
\eta < -1$.  Then $\psi > \theta_1^2$ and $P_i(\psi) > 0$ for $0 \le i \le D$.  
\item[{\rm (ii)}]  Assume $-1 < \eta < {\tilde \theta}_d$.  Then $ \psi < \theta_{d}^{2}$ and
$(-1)^{\lfloor \frac{i}{2} \rfloor}P_{i}(\psi) > 0$ for $0 \le i 
\le D$.
\item[{\rm (iii)}]  $P_i(\psi) \not= 0$ for $0 \leq i \leq D$.
\end{description}  
\end{lemma}
{\it Proof.}  (i)  Combining the inequalities ${\tilde \theta}_1 < 
\eta < -1$ with (\ref{tildeth1d}), 
(\ref{psidef}), and using Lemma \ref{bipeigs}, we routinely find $\psi > \theta_{1}^{2}$.  Thus 
$P_{i}(\psi) > 0$ $\;(0 \le i \le D)$ by Lemma \ref{Pnot0}(i).  \\
\noindent (ii)
Combining the inequalities $-1 < \eta < {\tilde \theta}_d$ with (\ref{tildeth1d}),
(\ref{psidef}), and using Lemma \ref{bipeigs}, we routinely find $ \psi < \theta_{d}^{2}$.  Thus 
$(-1)^{\lfloor \frac{i}{2} \rfloor} P_{i}(\psi) > 0$ $\;(0 \le i \le D)$ by Lemma \ref{Pnot0}(ii),(iii).  \\
\noindent (iii)
Immediate from (i), (ii) above.
\hfill $\Box$\\

\noindent Referring to Theorem \ref{M5}, 
we now consider a second basis for $Mv$.

\begin{definition} 
\label{M1}  \rm
With reference to Definition \ref{A}, let
$v$ denote a nonzero vector in $U$. Assume $v$ is an eigenvector
for $E^*_2A_{2}E^*_2$ and let $\eta$ denote the corresponding
eigenvalue. Assume 
$ {\tilde \theta}_1 <  
\eta < {\tilde \theta}_{d}$. 
We define the vectors $v_0, v_1, \ldots, v_{D-2}$ as follows:
\begin{description}
\item[{\rm (i)}]  Suppose $\eta \not=-1$. Then 
\begin{equation}\label{M1a}
v_i = \sum^{i}_{{h=0}\atop{ i-h { \mbox{ \tiny even}}}}
\frac{P_h(\psi)}{P_i(\psi)}\,
\frac{k_ib_ib_{i+1}}{k_hb_hb_{h+1}} \,p_h(A)v
\qquad \qquad (0 \leq i \leq D-2),
\end{equation}
where $\psi$ is from (\ref{psidef}).
\item[{\rm (ii)}]  Suppose $\eta =-1$. Then
$v_i = p_i(A)v$ for 
$0 \leq i \leq D-2$.
\end{description}
(The polynomials $p_i$ are from  
(\ref{PIPOLY}), and the $P_{i}$ are from (\ref{Pidef2}).)
\end{definition}

\begin{theorem} 
\label{M2}
With reference to Definition \ref{A}, let
$v$ denote a nonzero vector in $U$. Assume $v$ is an eigenvector
for $E^*_2A_{2}E^*_2$ and let $\eta$ denote the corresponding
eigenvalue. Assume
$ {\tilde \theta}_1 <  
\eta < {\tilde \theta}_d$. 
Then the vectors 
$v_0, v_1, \ldots , v_{D-2}$
from Definition \ref{M1} form a basis for $Mv$.
\end{theorem}
{\it Proof.} 
By  Theorem
\ref{M5}
  we find $Mv$ has dimension $D-1$. By this
and since $A$ generates $M$, we find $Mv$ has a  basis
$v, Av, \ldots, A^{D-2}v$.
For $0 \leq i \leq D-2$ the vector  $v_i$ is contained in the span of
$v, Av, \ldots,  A^iv$ but not in the span of
$v, Av, \ldots,  A^{i-1}v$.
It follows that $v_0, v_1, \ldots, v_{D-2}$ form a basis for $Mv$.
\hfill $\Box$\\

\noindent With reference to
Definition
\ref{M1},
we will show that the vectors 
$v_0, v_1, \ldots , v_{D-2}$
are mutually orthogonal and we will compute their 
square-norms. To do this we need the following result. 

\begin{theorem} 
\label{M7}
With reference to Definition \ref{A}, let
$v$ denote a nonzero vector in $U$. Assume $v$ is an eigenvector
for $E^*_2A_{2}E^*_2$ and let $\eta$ denote the corresponding
eigenvalue. Assume 
$ {\tilde \theta}_1 <  
\eta < {\tilde \theta}_d$. 
Let the vectors $v_0, v_1, \ldots, v_{D-2}$
be as in Definition
\ref{M1}.
\begin{description}
\item[{\rm (i)}]  Suppose $\eta \not=-1$. Then for $0 \leq i \leq D-2$ we have 
\begin{equation}\label{M7a}
v_i =   \sum_{j=1}^{D-1} g_i(\theta_j) E_jv,
\end{equation}
where
\begin{equation}\label{M7b}
g_i = \sum_{{h=0}\atop{ i-h { \mbox{ \tiny even}}}}^{i} 
\frac{P_h(\psi)}{P_i(\psi)}\,
\frac{k_ib_ib_{i+1}}{k_hb_hb_{h+1}} \,p_h
\end{equation}
and $\psi$ is from (\ref{psidef}).
\item[{\rm (ii)}]  Suppose $\eta =-1$. Then
$$
v_i =   \sum_{j=1}^{D-1} p_i(\theta_j) E_jv
\qquad \qquad (0 \leq i \leq D-2).
$$
\end{description}
\end{theorem}
{\it Proof.} (i) Let the integer $i$ be given. Comparing 
(\ref{M1a}),
(\ref{M7b}) we find $v_i=g_i(A)v$. Using this and 
(eii) we routinely obtain
$v_i =
\sum_{j=0}^D g_i(\theta_j) E_jv$.
Line  (\ref{M7a}) follows since $E_0v=0$, $E_Dv=0$ by Lemma
\ref{EJV}(i). \\
\noindent (ii) Similar to the proof of (i) above. 
\hfill $\Box$\\

\begin{theorem} 
\label{M3}
With reference to Definition \ref{A}, let
$v$ denote a nonzero vector in $U$. Assume $v$ is an eigenvector
for $E^*_2A_{2}E^*_2$ and let $\eta$ denote the corresponding
eigenvalue. Assume 
$ {\tilde \theta}_1 <  
\eta < {\tilde \theta}_d$. 
Then the vectors
 $v_0, v_1, \ldots, v_{D-2}$
from Definition
\ref{M1}
 are mutually orthogonal.
Moreover the square-norms of these vectors are given
as follows:
\begin{description}
\item[{\rm (i)}]  Suppose $\eta \not=-1$. Then 
\begin{equation}\label{M3a}
\|v_i\|^2 =  \frac{k_{i}b_{i}b_{i+1}c_{i+1}c_{i+2}}{k b_{1}(\psi - b_{2})} \,
\frac{P_{i+2}(\psi)}
{P_i(\psi)} 
\, \|v \|^2
\qquad \qquad (0 \leq i \leq D-2),
\end{equation}
where $\psi$ is from (\ref{psidef}).
\item[{\rm (ii)}]  Suppose $\eta =-1$. Then
$$
\|v_i\|^2 =  \frac{k_{i} b_i b_{i+1}}{k b_{1}} \,\|v\|^2
\qquad \qquad (0 \leq i \leq D-2).
$$
\end{description}
\end{theorem}
{\it Proof.} (i)
Let the 
polynomials $g_0, g_1, \ldots, g_{D-2}$ be as in
(\ref{M7b}). Using in order
Theorem 
\ref{M7},
Theorem \ref{M6}, 
and 
Theorem \ref{P10},
we find that for $0 \leq i,j\leq D-2$, 
\begin{eqnarray}
\langle v_i, v_j \rangle &=&  
\sum_{h=1}^{D-1} g_i(\theta_h)g_j(\theta_h)\|E_hv\|^2
\nonumber \\
&=& \sum_{h=1}^{D-1} g_i(\theta_h)g_j(\theta_h)\,
\frac{m_h(\theta_h 
-k)(\theta_{h}+k)(\theta_h^{2}-\psi)}{|X|kb_{1}(\psi -b_{2})} \|v\|^2
\nonumber \\
&=& \delta_{ij} \, \frac{k_{i}b_{i}b_{i+1}c_{i+1}c_{i+2}}{kb_{1}(\psi - b_{2})} \,
\frac{P_{i+2}(\psi)}
{P_i(\psi)}
\, \|v \|^2.
\nonumber
\end{eqnarray}
Apparently $v_0, v_1, \ldots, v_{D-2}$ are mutually orthogonal
and satisfy 
(\ref{M3a}). 
\\

\noindent (ii) 
The argument is similar  to (i) above,
with the $p_i$ taking the place of the $g_i$
and 
Lemma  \ref{PIPJ} taking the place of
Theorem \ref{P10}.
\hfill $\Box$\\

\begin{theorem} 
\label{M4}
With reference to Definition \ref{A}, let
$v$ denote a nonzero vector in $U$. Assume $v$ is an eigenvector
for $E^*_2A_{2}E^*_2$ and let $\eta$ denote the corresponding
eigenvalue. Assume 
$ {\tilde \theta}_1 <  
\eta < {\tilde \theta}_d$. 
With respect to the basis 
$v_0, v_1, \ldots, v_{D-2}$
for $Mv$ given in Definition
\ref{M1}
the matrix representing $A$ is
$$
\left(\begin{array}{cccccc}
0 & \omega_1 & & & & {\bf 0}\\
c_1 & 0 & \omega_2 & & & \\
 & c_2 & \cdot & \cdot & & \\
 & & \cdot & \cdot &  \cdot& \\
 & & & \cdot & \cdot & \omega_{D-2} \\
 {\bf 0} & & & & c_{D-2} & 0
\end{array} \right),
$$
where the $\omega_i$ are as follows:
\begin{description}
\item[{\rm (i)}]  Suppose $\eta \not=-1$. Then 
\begin{equation}
\omega_i = 
\label{M4aa}  
\,\frac{ b_{i+1}c_{i+2}}{c_i}
\frac{P_{i-1}(\psi)P_{i+2}(\psi)}
{P_{i}(\psi) P_{i+1}(\psi)}
\qquad \qquad (1 \leq i \leq D-2),
\end{equation}
where $\psi$ is from (\ref{psidef}).
\item[{\rm (ii)}]  Suppose $\eta =-1$. Then
\begin{equation}
\omega_i = b_{i+1} \qquad \qquad (1 \leq i \leq D-2).
\nonumber
\end{equation}
\end{description}
\end{theorem}
{\it Proof.} (i) For  $0 \leq i \leq D-2$ we define $g_i$ as in
(\ref{M7b}). 
 Setting 
$\lambda =A$ and $\theta =\psi$ in Theorem
\ref{P9} we find
\begin{equation}\label{M4d}
Ag_i(A) = c_{i+1}g_{i+1}(A) +  \omega_ig_{i-1}(A) 
\qquad \qquad (0 \le i \le D-2),
\end{equation}
where $g_{-1}=0$, $\omega_{0}=0$, $g_{D-1}=p_{D-1}$, and the $\omega_i$ are
from 
(\ref{M4aa}).
 Observe $g_i(A)v=v_i$ for $0 \leq i \leq D-2$. 
 Applying both equations in (\ref{PDAeqs}) to $v$ and recalling $Jv=0$, 
 $J'v=0$, we find $p_{D-1}(A)v=0$.
Applying (\ref{M4d}) to $v$, and simplifying the result using
these comments, we find
$$
Av_i = c_{i+1}v_{i+1} + \omega_iv_{i-1} 
\qquad \qquad (0 \le i \le D-2),
$$
where $v_{-1}=0$ and $v_{D-1}=0$. 
The result follows.

\noindent (ii) The argument is similar to (i) above, with the
$p_i$ taking the place of the $g_i$ and 
(\ref{PIRECUR1}) taking the place of
Theorem \ref{P9}.
\hfill $\Box$\\


\section{The thin irreducible $T$-modules with endpoint 2
and local eigenvalue $\eta$ 
$( {\tilde \theta}_1 <  
\eta < {\tilde \theta}_d$)}  \label{main1}

\noindent 
With reference to Definition \ref{A}, we now
describe the thin irreducible $T$-modules
with endpoint 2
and local eigenvalue $\eta$
$( {\tilde \theta}_1 <  
\eta < {\tilde \theta}_d$). 
This section contains some of our
main results. Because of this 
we have tried to make it as self-contained as 
possible.

\begin{theorem} 
\label{T5}
With reference to Definition \ref{A}, let
$W$ denote a thin irreducible $T$-module with  endpoint 2
and local eigenvalue $\eta$
$( {\tilde \theta}_1 <  
\eta < {\tilde \theta}_d)$. 
Let $v$ denote
a nonzero vector in $E^*_2W$. 
Then $W=Mv$. The vectors  
\begin{equation}\label{T5a}
E_1v, E_2v, \ldots, E_{D-1}v 
\end{equation}
form a basis for $W$ and $E_0v=0$, $E_{D}v=0$.
\end{theorem}
{\it Proof.} To see $W=Mv$, observe that $W$ contains $v$ and is
invariant under $M$ so $Mv \subseteq W$. We assume $W$ is thin
with endpoint 2, so the dimension of $W$ is at most $D-1$.
By Definition 
\ref{TM}
the vector $v$ is contained in $U$. Moreover $v$ is an eigenvector
for $E^*_2A_{2}E^*_2$ with eigenvalue $\eta$.
Now Theorem \ref{M5} applies. By that theorem
 $Mv$ has dimension $D-1$ so $W=Mv$. The remaining
assertions of the present theorem follow in view of 
Theorem \ref{M5}.
\hfill $\Box$\\

\begin{theorem} 
\label{T6}
With reference to Definition \ref{A}, let
$W$ denote a thin irreducible $T$-module with  endpoint 2
and local eigenvalue $\eta$
$( {\tilde \theta}_1 <  
\eta < {\tilde \theta}_d)$. 
Then the basis vectors for $W$  from   
(\ref{T5a})  
are mutually orthogonal.
Moreover the square-norms of these vectors are given
as follows:
\begin{description}
\item[{\rm (i)}]  Suppose $\eta \not=-1$. Then 
$$\|E_iv\|^2
= \frac{m_i(\theta_i-k)(\theta_{i}+k)(\theta_i^{2}-\psi)}
{|X|k b_{1}(\psi - b_{2})}
\|v\|^2
\qquad \qquad (1 \le i \le D-1),
$$ where \begin{equation}  \label{psidef2}  \psi = b_{2}\left(1- 
\frac{b_{3}}{1+\eta}\right). \end{equation}
\item[{\rm (ii)}]  Suppose $\eta =-1$. Then
$$
\|E_iv\|^2
= \frac{m_i(k - \theta_i)(k + \theta_i)}{|X|kb_1} \|v\|^2
\qquad \qquad (1 \le i \le D-1).
$$
\end{description}
(The scalar $m_i$ denotes the multiplicity of $\theta_i$.)
\end{theorem}
{\it Proof.}  
By Definition \ref{TM} the vector $v$ is contained in $U$. 
Moreover $v$ is an eigenvector
for $E^*_2A_{2}E^*_2$ with eigenvalue $\eta$.
Applying  Theorem \ref{M6} we obtain the result.
\hfill $\Box$\\

\begin{theorem} 
\label{T1}
With reference to Definition \ref{A}, let
$W$ denote a thin irreducible $T$-module with  endpoint 2
 and local eigenvalue $\eta$
$( {\tilde \theta}_1 <  
\eta < {\tilde \theta}_d)$. 
 Let $v$ denote
a nonzero vector in $E^*_2W$.
\begin{description}
\item[{\rm (i)}]  Suppose $\eta \not=-1$. Then 
$$
E^*_{i+2}A_iv = \sum_{{h=0}\atop{ i-h { \mbox{ \tiny even}}}}^{i} 
\frac{P_h(\psi)}{P_i(\psi)}\,
\frac{k_ib_ib_{i+1}}{k_hb_hb_{h+1}} \,p_h(A)v
\qquad \qquad (0 \leq i \leq D-2),
$$
where $\psi$ is from (\ref{psidef2}).
\item[{\rm (ii)}]  Suppose $\eta =-1$. Then
$$
E^*_{i+2}A_iv = p_i(A)v
\qquad \qquad (0 \leq i \leq D-2).
$$
\end{description}
(The polynomials $p_i$ are from 
(\ref{PIPOLY}), and the $P_{i}$ are from (\ref{Pidef2}).)
\end{theorem}
{\it Proof.}
By Definition \ref{TM} the vector $v$ is contained in $U$. 
Moreover $v$ is an eigenvector
for $E^*_2A_{2}E^*_2$ with eigenvalue $\eta$.
Let the vectors $v_0, v_1, \ldots, v_{D-2}$
be as in Definition 
\ref{M1}. 
We show $E^*_{i+2}A_iv=v_i$ for $0 \leq i \leq D-2$.
Using 
(\ref{REL3})
we find 
$A^iv$ is contained in $E^*_2W+\cdots + E^*_{i+2}W$
for $0 \leq i \leq D-2$. Also for 
 $0 \leq i \leq D-2$, 
 $v_i$ is a linear combination
of $v, Av, \ldots, A^iv $, so $v_i $ is contained in 
$E^*_2W+\cdots + E^*_{i+2}W$.
By this and since  $v_0, v_1, \ldots, v_{D-2}$
are  linearly independent, we find 
\begin{equation}\label{mp}
v_0, v_1, \ldots, v_i \quad \hbox{is a basis for }   \quad 
E^*_2W+E^*_3W+\cdots + E^*_{i+2}W \qquad  (0 \leq i \leq D-2).
\end{equation}
For the rest of this proof, 
fix an integer $i$ $(0 \leq i \leq D-2)$. 
We show that $v_i $ is contained in $E^*_{i+2}W$.
To see this,
recall $E^*_2W,\ldots, E^*_DW$
are mutually orthogonal. Therefore $E^*_{i+2}W $ is equal to 
the orthogonal complement of 
$E^*_2W+\cdots + E^*_{i+1}W $ in 
$E^*_2W+\cdots + E^*_{i+2}W$. Recall $v_i$ is orthogonal to 
each of $v_0, v_1, \ldots, v_{i-1}$. By (\ref{mp}) the
vectors 
 $v_0, v_1, \ldots, v_{i-1}$ form a basis for  
$E^*_2W+\cdots + E^*_{i+1}W $ so $v_i$ is orthogonal to  
$E^*_2W+\cdots + E^*_{i+1}W $. Apparently $v_i$ is 
contained in $E^*_{i+2}W$ as desired.
We show that $E^*_{i+2}A_iv=v_i$.
We mentioned that
the vector $v_i$ is a linear combination 
of $v, Av, \ldots, A^iv$. In this combination 
the coefficient of $A^iv$ is $(c_1 c_2\cdots c_i)^{-1}$ 
in view of Lemma \ref{P8}(ii). Similarly $A_iv$ 
is a linear combination of
$v, Av, \ldots, A^iv$, and in this combination   
the coefficient of $A^iv$ is $(c_1 c_2\cdots c_i)^{-1}$. 
Apparently $A_iv-v_i$ is a linear combination of 
$v, Av, \ldots, A^{i-1}v$. From this and our above comments
$A_iv-v_i$ is contained in  
$E^*_2W+\cdots + E^*_{i+1}W $ so $E^*_{i+2}(A_iv-v_i)$ is
zero. We already  showed that $v_i \in E^*_{i+2}W$ so
 $E^*_{i+2}v_i=v_i$. Now
$E^*_{i+2}A_iv=v_i$ as desired.   
\hfill $\Box$\\

\begin{theorem} 
\label{T2}
With reference to Definition \ref{A}, let
$W$ denote a thin irreducible $T$-module with  endpoint 2
 and local eigenvalue $\eta$
$( {\tilde \theta}_1 <  
\eta < {\tilde \theta}_d)$. 
 Let $v$ denote
a nonzero vector in $E^*_2W$. 
Then the vectors  
\begin{equation}\label{T2a}
E^*_{i+2}A_iv \qquad \qquad (0 \leq i \leq D-2)
\end{equation}
form a basis for $W$.
\end{theorem}
{\it Proof.}  By Definition 
\ref{TM} the vector $v$ is contained in $U$. 
Moreover $v$ is an eigenvector for $E^*_2A_{2}E^*_2$
with eigenvalue $\eta$. Let the vectors $v_0, v_1, \ldots, v_{D-2}$
be as in Definition
\ref{M1}.  By 
 Theorem
\ref{M2} the vectors 
$v_0, v_1, \ldots, v_{D-2}$
form a basis for $Mv$. Recall $Mv=W$ by Theorem
\ref{T5} so 
$v_0, v_1, \ldots, v_{D-2}$
form a basis for $W$. By
Theorem \ref{T1}
$v_i = E^*_{i+2}A_iv$ for $0 \leq i \leq D-2$ and
the result follows.
\hfill $\Box$\\

\begin{theorem} 
\label{T3}
With reference to Definition \ref{A}, let
$W$ denote a thin irreducible $T$-module with  endpoint 2
and local eigenvalue $\eta$
$( {\tilde \theta}_1 <  
\eta < {\tilde \theta}_d)$. 
Then the vectors in (\ref{T2a}) are mutually  orthogonal.
Moreover the square-norms of these vectors are given
as follows:
\begin{description}
\item[{\rm (i)}]  Suppose $\eta \not=-1$. Then 
$$
\|E^*_{i+2}A_iv\|^2 =   \frac{k_{i}b_{i}b_{i+1}c_{i+1}c_{i+2}}
{k b_{1}(\psi - b_{2})} \,
\frac{P_{i+2}(\psi)}
{P_i(\psi)} 
\, \|v \|^2
\qquad \qquad (0 \leq i \leq D-2),
$$
where $\psi$ is from (\ref{psidef2}).
\item[{\rm (ii)}]  Suppose $\eta =-1$. Then
$$
\|E^*_{i+2}A_iv\|^2 =  \frac{k_{i} b_i b_{i+1}}{k b_{1}} \,\|v\|^2
\qquad \qquad (0 \leq i \leq D-2).
$$
\end{description}
\end{theorem}
{\it Proof.}
  By Definition 
\ref{TM} the vector $v$ is contained in $U$. 
Moreover $v$ is an eigenvector for $E^*_2A_{2}E^*_2$
with eigenvalue $\eta$. 
The result follows in view of 
Theorem  
\ref{M3} and Theorem
\ref{T1}.
\hfill $\Box$\\

\begin{theorem} 
\label{T4}
With reference to Definition \ref{A}, let
$W$ denote a thin irreducible $T$-module with  endpoint 2
 and local eigenvalue $\eta$
$( {\tilde \theta}_1 <  
\eta < {\tilde \theta}_D)$. 
With respect to the basis for $W$ given in 
(\ref{T2a}) the matrix representing $A$ is
$$
\left(\begin{array}{cccccc}
0 & \omega_1 & & & & {\bf 0}\\
c_1 & 0 & \omega_2 & & & \\
 & c_2 & \cdot & \cdot & & \\
 & & \cdot & \cdot &  \cdot& \\
 & & & \cdot & \cdot & \omega_{D-2} \\
 {\bf 0} & & & & c_{D-2} & 0
\end{array} \right),
$$
where the $\omega_i$ are as follows. 
\begin{description}
\item[{\rm (i)}]  Suppose $\eta \not=-1$. Then 
\begin{equation} \nonumber
\omega_i =  
\,\frac{ b_{i+1}c_{i+2}}{c_i}
\frac{P_{i-1}(\psi)P_{i+2}(\psi)}
{P_{i}(\psi) P_{i+1}(\psi)}
\qquad \qquad (1 \leq i \leq D-2),
\end{equation}
where $\psi$ is from (\ref{psidef2}).
\item[{\rm (ii)}]  Suppose $\eta =-1$. Then
\begin{equation}
\omega_i = b_{i+1} \qquad \qquad (1 \leq i \leq D-2).
\nonumber
\end{equation}
\end{description}
\end{theorem}
{\it Proof.} 
  By Definition 
\ref{TM} the vector $v$ is contained in $U$. 
Moreover $v$ is an eigenvector for $E^*_2A_{2}E^*_2$
with eigenvalue $\eta$. 
The result follows in view of 
Theorem  
\ref{M4} 
and Theorem
\ref{T1}.
\hfill $\Box$\\

\begin{theorem} 
\label{T7}
With reference to Definition \ref{A}, let
$W$ denote a thin irreducible $T$-module with  endpoint 2
 and local eigenvalue $\eta$
$( {\tilde \theta}_1 <  
\eta < {\tilde \theta}_d)$. 
Let $v$ denote a nonzero vector in $E^*_2W$.
\begin{description}
\item[{\rm (i)}]  Suppose $\eta \not=-1$. Then for $0 \leq i \leq D-2$ we have 
$$
E^*_{i+2}A_iv =   \sum_{j=1}^{D-1} g_i(\theta_j) E_jv,
$$
where
$$
g_i = \sum_{{h=0}\atop{ i-h { \mbox{ \tiny even}}}}^{i} 
\frac{P_h(\psi)}{P_i(\psi)}\,
\frac{k_ib_ib_{i+1}}{k_hb_hb_{h+1}} \,p_h
$$
and $\psi$ is from (\ref{psidef2}).
\item[{\rm (ii)}]  Suppose $\eta =-1$. Then
$$
E^*_{i+2}A_iv =   \sum_{j=1}^{D-1} p_i(\theta_j) E_jv
\qquad \qquad (0 \leq i \leq D-2).
$$
\end{description}
\end{theorem}
{\it Proof.} 
  By Definition 
\ref{TM} the vector $v$ is contained in $U$. 
Moreover $v$ is an eigenvector for $E^*_2A_{2}E^*_2$
with eigenvalue $\eta$. 
The result follows in view of 
 Theorem \ref{M7} and  
Theorem \ref{T1}.
\hfill $\Box$\\

\noindent In summary we have the following theorem.

\begin{theorem}
\label{SUM}
With reference to Definition \ref{A}, let
$W$ denote a thin irreducible $T$-module with  endpoint 2
 and local eigenvalue $\eta$
$( {\tilde \theta}_1 <  
\eta < {\tilde \theta}_d)$. 
Then $W$ has dimension $D-1$.
For $0 \leq i \leq D$, $E^*_iW$ is zero if $i \in \{0,1 \}$ and has dimension
1 if $2 \leq i \leq D$. Moreover  
 $E_iW$ is zero if $i \in \{0, D \}$ and has dimension
1 if $1 \leq i \leq D-1$. 
\end{theorem}
{\it Proof.} The dimension of $W$ is 
$D-1$ by Theorem
\ref{T5}.
Fix an integer $i$ $(0 \leq i \leq D)$.
From Theorem
\ref{T2} we find 
 $E^*_iW$ is zero if $i \in \{0, 1 \}$ and has dimension
1 if $2 \leq i \leq D$.   
From Theorem
\ref{T5}
we find 
 $E_iW$ is zero if $i \in \{0, D \}$ and has dimension
1 if $1 \leq i \leq D-1$.   
\hfill $\Box$\\


\section{Some multiplicities}

With reference to Definition \ref{A}, let $W$ denote a thin
irreducible $T$-module with endpoint 2 and local eigenvalue $\eta$.
In this section we show that the isomorphism class of $W$ as a $T$-module
is determined
by $\eta$. We show that the 
 multiplicity with which $W$ appears 
in the standard module $V$ is at most the number of times 
$\eta$ appears among $\eta_{k+1}, \eta_{k+2}, \ldots, \eta_{k_{2}}$.
We investigate the case of equality.

\begin{theorem} 
\label{T8}
With reference to Definition \ref{A}, let
$W$ denote a thin irreducible $T$-module with  endpoint 2
and local eigenvalue $\eta$.
Let $W'$ denote an irreducible $T$-module. Then the following
(i), (ii) are equivalent:
\begin{description}
\item[{\rm (i)}] $W $ and $W'$ are isomorphic as $T$-modules.
\item[{\rm (ii)}]
$W'$ is thin  with  endpoint 2 
and local eigenvalue $\eta$.
\end{description}
\end{theorem}
{\it Proof.} (i)$\Rightarrow $(ii) Clear.

\noindent (ii)$\Rightarrow $(i)  First observe that $\eta$ satisfies one of the cases (i)--(iv) mentioned below Definition \ref{TM}.  If $\eta$ satisfies case (i) or case (ii) then statement (i) of the present theorem holds by 
\cite[Theorem 14.1]{maclean5}.  Now assume $\eta$ satisfies case (iii) or case (iv).  For notational convenience set $e=1$ if $\eta$ satisfies case (iii) and set $e=0$ if $\eta$ satisfies case (iv).  We display an isomorphism of $T$-modules 
from $W$ to $W^\prime $. 
Observe $E^*_2W$ and $E^*_2W'$ are both nonzero.
Let $v$ (resp. $v'$) denote a nonzero
vector in $E^*_2W$ (resp. $E^*_2W'$).  
By Theorem
\ref{T2d}  or Theorem \ref{T2} the vectors
\begin{equation}
E^*_{i+2}A_iv \qquad \qquad  (0 \leq i \leq D-2-e)
\label{b1}
\end{equation} 
form
a basis for $W$. Similarly the vectors
\begin{equation}
E^*_{i+2}A_iv' \qquad \qquad  (0 \leq i \leq D-2-e)
\label{b1p}
\end{equation} 
form
a basis for $W'$.
Let $\sigma : W \rightarrow W'$ denote the isomorphism of
vector spaces that sends 
 $E^*_{i+2}A_iv$ 
 to $E^*_{i+2}A_iv'$ for
$0 \leq i \leq D-2-e$. We show $\sigma $ is an isomorphism
of $T$-modules. By Theorem
\ref{T4d} or Theorem \ref{T4}
 the matrix representing $A$ with respect to the basis
(\ref{b1}) 
  is equal to the matrix 
 representing $A$ with respect to the basis
  (\ref{b1p}).
  It follows $ \sigma A - A \sigma $ vanishes on $W$. 
  From  the construction we find
  that  for $0 \leq h \leq D$,
 the matrix representing $E^*_h$ with respect to the basis
 (\ref{b1}) 
is equal to the matrix 
 representing $E^*_h$ with respect to the basis
 (\ref{b1p}). 
  It follows $ \sigma E^*_h  -  E^*_h \sigma $ vanishes on $W$. 
 The algebra $T$ is generated by 
  $A,
E_0^*,E_1^*,\ldots,E_D^*$.
 It follows  $\sigma B - B \sigma $ vanishes
on $W$ for all $B \in T$. We now see $\sigma $ is an isomorphism
of $T$-modules from $W$ to $W^\prime $.
\hfill $\Box$\\

\begin{lemma}
\label{T11}
With reference to Definition \ref{A}, for all $\eta \in \R$ we have  
\begin{equation} \label{ueta}
U_\eta \supseteq E^*_2H_\eta,
\end{equation}
where $H_\eta $ denotes the subspace of $V$ spanned by 
all the thin 
irreducible $T$-modules with endpoint 2
and local eigenvalue $\eta$.
\end{lemma}
{\it Proof.}
Observe $E^*_2H_\eta $ is spanned by the $E^*_2W$, where $W$ ranges
over all the thin irreducible $T$-modules with endpoint 2 and local
eigenvalue $\eta$. For all such $W$ the space $E^*_2W$
is contained in $U_\eta $ by 
Definition
\ref{TM}. The result follows.
\hfill $\Box$\\

\noindent We remark on the dimension of the right-hand side in 
(\ref{ueta}). To do this we make a definition.

\begin{definition}
\label{MULTD}  \rm
With reference to Definition \ref{A}, and from our discussion
in Section \ref{subconsection}, 
 the standard module $V$ can be
decomposed into an orthogonal direct sum of irreducible $T$-modules.
Let $W$ denote an irreducible $T$-module. By the {\it multiplicity
with which $W$ appears in $V$}, we mean the number of irreducible
$T$-modules in the above decomposition which are isomorphic to $W$. 
\end{definition}

\begin{definition}
\label{D11A}  \rm
With reference to Definition \ref{A}, for all $\eta \in \R$
we let $\mu_\eta$ denote the multiplicity with which
$W$ appears in $V$, where $W$ is a thin irreducible $T$-module
with endpoint 2 and local eigenvalue $\eta$.
If no such $W$ exists we interpret $\mu_\eta =0$.
\end{definition}

\begin{theorem} 
\label{T12}
With reference to Definition \ref{A}, for all $\eta \in \R$ the
following scalars are equal:
\begin{description}
\item[{\rm (i)}]
The scalar $\mu_\eta $ 
from Definition
\ref{D11A}.
\item[{\rm (ii)}]
The dimension of $E^*_2H_\eta $,
where $H_\eta $ is from
Lemma \ref{T11}.
\end{description}
Moreover
\begin{equation}
 \hbox{mult}_\eta  \geq \mu_\eta.
\label{43new}
\end{equation}
\end{theorem}
{\it Proof.} We first show that
 $\mu_\eta $  is equal to 
 the dimension of 
$E^*_2H_\eta $.
Observe $H_\eta$ is a $T$-module so it is an orthogonal direct sum
of irreducible $T$-modules. More precisely 
\begin{equation}\label{hometa}
H_\eta = W_1 + W_2 +\cdots + W_m \qquad \qquad {\rm (orthogonal\ direct\ sum)},
\end{equation}
where $m$ is a nonnegative integer, and 
 where $W_1,W_2,\ldots,W_m$ are thin irreducible $T$-modules 
with endpoint 2 and local eigenvalue $\eta$.
Apparently $m$ is equal to  $\mu_\eta$.
We show $m$ is equal to the dimension of $E^*_2H_\eta $.
Applying $E_2^*$ to (\ref{hometa}) 
 we find
\begin{equation}\label{hometa2}
E^*_2H_\eta  = E_2^*W_1 + E_2^*W_2 +\cdots + E_2^*W_m \qquad \qquad {\rm (orthogonal\ direct\ sum)}.
\end{equation}
 Observe each  
summand on the right in (\ref{hometa2}) has dimension 1.
These summands are mutually orthogonal so $m$ is equal to the dimension of
 $E^*_2H_\eta $.
 Now  
$\mu_\eta $ is equal to the dimension of $E^*_2H_\eta $.
We mentioned earlier that
the dimension of $U_\eta$ is $\hbox{mult}_\eta$.
Combining these facts  with
Lemma \ref{T11} we obtain
(\ref{43new}).
\hfill $\Box$\\

\noindent 
We are interested in the case of equality in 
(\ref{ueta}) and
(\ref{43new}).
 We begin with a result which is a routine consequence
 of 
Lemma \ref{lem:mvisshort}.

\begin{lemma} \cite[Lemma 14.2]{maclean5}
With reference to Definition \ref{A}, choose
$n \in \{1,d\}$ if $D$ is odd, and let $n=1$ if $D$ is even.
Let $\eta ={\tilde \theta}_n$. Then $U_\eta  = E^*_2H_\eta$ 
and 
$\hbox{mult}_\eta = \mu_\eta$.
\end{lemma}

\begin{lemma} 
\label{N}
With reference to Definition \ref{A},
let $L$ denote the subspace of $V$ spanned by the nonthin
irreducible $T$-modules
with endpoint 2. Then
\begin{equation}
U = E^*_2L + \sum_{\eta \in \Phi} E^*_2H_\eta \qquad \qquad
(\hbox{orthogonal direct sum}).
\label{Neq}
\end{equation}
\end{lemma}
{\it Proof.} Let $S$ denote the subspace of
$V$ spanned by all 
irreducible $T$-modules
with endpoint 2, thin or not. Then
\begin{equation}
\label{Neq2}
S = L + \sum_{\eta \in \Phi}H_\eta \qquad \qquad
(\hbox{orthogonal direct sum}).
\end{equation}
Applying $E^*_2$ to each term in  
(\ref{Neq2}) and using $E^*_2S=U$ we obtain
(\ref{Neq}).
\hfill $\Box$\\

\begin{theorem}
\label{T13}
With reference to Definition \ref{A}, 
the following (i)--(iii) are equivalent:
\begin{description}
\item[{\rm (i)}] Equality holds in 
(\ref{ueta}) for all $\eta \in \R$.
\item[{\rm (ii)}] Equality holds in 
(\ref{43new}) for all $\eta \in \R$.
\item[{\rm (iii)}] Every irreducible $T$-module with endpoint 2 is thin. 
\end{description}
\end{theorem}
{\it Proof.} 
\noindent 
(i)$\Leftrightarrow $(ii) Recall $\hbox{mult}_\eta $ (resp. 
$\mu_\eta $) is the dimension of
 $U_\eta $ (resp. 
$E^*_2H_\eta$).
\\
\noindent
(i)$\Rightarrow $(iii) 
Let $W$ denote an irreducible $T$-module with endpoint 2.
We show $W$ is thin. Suppose not. Then
$W$ is contained in 
the space $L$ from 
Lemma 
\ref{N}.  Observe $E^*_2W\not=0$ since $W$ has endpoint 2,
so $E^*_2L\not=0$.
We show 
$E^*_2L=0$ to get a contradiction.
We assume 
$
U_\eta
=
E^*_2H_\eta
$
for all $\eta \in \R$;
combining this with 
(\ref{Ubreakdown}) we find 
$U=\sum_{\eta \in \Phi}E^*_2H_\eta $.
From this and 
Lemma \ref{N} we find $E^*_2L=0$. We 
now have a contradiction and it follows $W$ is thin. \\
\noindent
(iii)$\Rightarrow $(i) 
There does not exist a nonthin irreducible $T$-module with endpoint 2, so
$L=0$. Setting $L=0$ in 
(\ref{Neq})
we find
$U=\sum_{\eta \in \Phi}E^*_2H_\eta $. Combining this with 
(\ref{Ubreakdown}) and 
Lemma
\ref{T11} we routinely find
$
U_\eta
=
E^*_2H_\eta
$
for all $\eta \in \Phi$.
For any real number $\eta $ that is not in $\Phi$
the spaces $U_\eta $ and $H_\eta$ are both 0. 
Now 
$
U_\eta
=
E^*_2H_\eta
$
for all $\eta \in \R$.
\hfill $\Box$\\

\bigskip

\noindent
Mark S. MacLean \\
Mathematics Department \\
Seattle University \\
901 Twelfth Avenue \\
Seattle WA  98122-1090 USA \\
email:  macleanm@seattleu.edu

\bigskip

\noindent
Paul Terwilliger \\
Mathematics Department \\
University of Wisconsin \\
480 Lincoln Drive \\
Madison WI  53706-1388 USA \\
email:  terwilli@math.wisc.edu


\begin{thebibliography}{35}


\bibitem{bannai} E. Bannai and T. Ito. {\it Algebraic Combinatorics I:
      Association Schemes.} Benjamin/Cummings, London, 1984.

\bibitem{biggs} N. Biggs. {\it Algebraic graph theory.} 
      Cambridge U. Press, London, 1994.


\bibitem{bcn} A. E. Brouwer, A. M. Cohen, and A. Neumaier. 
      {\it Distance-Regular Graphs.} Springer-Verlag, Berlin, 1989.




\bibitem{caugh2} J. S. Caughman IV. 
{\it The Terwilliger algebras of bipartite $P$-and $Q$-polynomial
association
schemes}.
Discrete Math. {\bf 196} (1999), 65--95. 

\bibitem{caugh3} J. S. Caughman IV and N. Wolff.  {\it The Terwilliger algebra of a distance-regular graph that supports a spin model}.  J. Alg. Combin.  {\bf 21}  (2005), 289--310.

\bibitem{caugh4} J. S. Caughman IV, M. S. MacLean, and P. Terwilliger.  {\it 
The Terwilliger algebra of an almost-bipartite $P$- and $Q$-polynomial association scheme}.
Discrete Math.  {\bf 292}  (2005), 17--44.


\bibitem{collins1} B. Collins. 
{\it The girth of a thin distance-regular graph}.
Graphs Combin. {\bf 13} (1997), 21--30. 

\bibitem{collins2} B. Collins.
{\it The Terwilliger algebra of an almost-bipartite distance-regular 
graph and its antipodal 2-cover}.
Discrete Math. {\bf 216} (2000), 35--69.

\bibitem{curtin1} B. Curtin.
{\it Bipartite distance-regular graphs I}.
Graphs Combin. {\bf 15} (1999), 
143--158.

\bibitem{curtin2} B. Curtin. 
{\it Bipartite distance-regular graphs II}.
Graphs Combin. {\bf 15} (1999), 
377--391.

\bibitem{curtin3} B. Curtin.
{\it 2-homogeneous bipartite distance-regular graphs}.
Discrete Math.
{\bf 187} (1998), 39--70.

\bibitem{curtin4} B. Curtin. 
{\it Distance-regular graphs which support a spin model are 
thin}. 16th British Combinatorial Conference (London, 1997). Discrete Math.
{\bf 197/198} (1999), 205--216.

\bibitem{curtin5} B. Curtin. 
{\it Almost 2-homogeneous bipartite distance-regular graphs}.
Europ. J. Combin.
{\bf 21} (2000), 865--876.

\bibitem{curtin9} B. Curtin.  {\it The Terwilliger algebra of a 2-homogeneous bipartite distance-regular graph}.  J. Combin. Theory Ser. B  {\bf 81}  (2001),  125--141.

\bibitem{curtin7} B. Curtin.  {\it Algebraic characterizations of graph regularity conditions}.  Des. Codes Cryptogr.  {\bf 3}4  (2005),  241--248.

\bibitem{curtin8} B. Curtin and K. Nomura.  {\it 1-homogeneous, pseudo-1-homogeneous, and 1-thin distance-regular graphs}.  J. Combin. Theory Ser. B  {\bf 93}  (2005),  279--302.



\bibitem{curtin6} B. Curtin and K. Nomura.
{\it Distance-regular graphs related to the quantum enveloping 
algebra of $sl(2)$}.
J. Alg. Combin.   
{\bf 12} (2000), 25--36.

\bibitem{CR} C. Curtis and I. Reiner.
{\it Representation theory of finite groups and associative algebras}.
Interscience, New York, 1962.

\bibitem{dickie1} G. Dickie. 
{\it Twice $\;Q$-polynomial distance-regular graphs are thin}.
Europ. J. Combin. {\bf 16} (1995), 555--560. 

\bibitem{dickie2} G. Dickie and P. Terwilliger.
{\it A note on thin $\;P$-polynomial and dual-thin
$\;Q$-polynomial symmetric association  schemes}.
J. Alg. Combin. {\bf 7} (1998), 5--15. 

\bibitem{egge1} E. Egge.
{\it A generalization of the Terwilliger algebra}. 
J. Algebra {\bf 233} (2000), 213--252.

\bibitem{egge2} E. Egge.
{\it The generalized Terwilliger algebra and its finite dimensional modules when $d=2$}.
J. Algebra {\bf 250} (2002), 178--216.

\bibitem{go} J. T.  Go.
{\it The Terwilliger algebra of the Hypercube $Q_D$}.
Europ. J. Combin. {\bf 23} (2002), 399-429.

\bibitem{go2} J. T. Go and P. Terwilliger. 
{\it Tight distance-regular graphs and the subconstituent algebra}.
Europ. J. Combin.  {\bf 23}  (2002),  793--816.

\bibitem{godsil} C. D. Godsil. {\it Algebraic Combinatorics.}
Chapman and Hall, Inc., New York, 1993. 

\bibitem{hobart} S. A. Hobart and T. Ito. {\it The structure of
nonthin irreducible $T$-modules: ladder bases and classical parameters}.
J. Alg. Combin. {\bf 7} (1998), 53--75.






\bibitem{jkt} A. Juri\v si\'c, J. Koolen and P. Terwilliger.
{\it Tight distance-regular graphs.} 
J. Alg. Combin. {\bf 12} (2000), 163--197. 
 


\bibitem{maclean5}
M. S. MacLean and P. Terwilliger. {\it Taut distance-regular graphs and 
the subconstituent 
algebra.} Discrete Math., accepted.











\bibitem{aap5} A. A. Pascasio.  {\it On the multiplicities of the primitive idempotents of a $Q$-polynomial distance-regular graph}.  Europ. J. Combin.  {\bf 23}  (2002),  1073--1078.

\bibitem{tanabe} K. Tanabe. {\it The irreducible modules of the Terwilliger
algebras of Doob schemes}. J. Alg. Combin. {\bf 6} (1997), 173--195.


\bibitem{terwSub1} P. Terwilliger. {\it The subconstituent algebra of
an association scheme I}. J. Alg. Combin. {\bf 1} (1992), 363--388.  

\bibitem{terwSub2} P. Terwilliger. {\it The subconstituent algebra of
an association scheme II}. J. Alg. Combin. {\bf 2} (1993), 73--103.  

\bibitem{terwSub3} P. Terwilliger. {\it The subconstituent algebra of
an association scheme III}. J. Alg. Combin. {\bf 3} (1993), 177--210.  

\bibitem{terw2} P. Terwilliger.  {\it An inequality involving the local eigenvalues of a distance-regular graph}.  J. Alg. Combin.  {\bf 19}  (2004), 143--172.

\bibitem{terw3} P. Terwilliger.  {\it The subconstituent algebra of a distance-regular graph; thin modules with endpoint one}.  Linear Algebra Appl.  {\bf 356}  (2002), 157--187.



\bibitem{tomiyama2} M. Tomiyama and N. Yamazaki. {\it The subconstituent
algebra of a strongly regular graph}. Kyushu J. Math. {\bf 48} (1998),
323--334.

\end{thebibliography}
\end{document}